\documentclass[12pt]{article}

\usepackage[T1,T2A]{fontenc}
\usepackage[utf8]{inputenc}

\usepackage{amsfonts}
\usepackage{amssymb}
\usepackage{amsmath}
\usepackage{amsthm}

\def \le {\leqslant}
\def \ge {\geqslant}

\usepackage{tikz}
\usetikzlibrary{calc,arrows,decorations.pathreplacing,fadings,3d,positioning}

\begin{document}
 
 \begin{Large}

 \centerline{\bf \"Uber die  Funktionen des Irrationalit\"atsma\ss es}
 \vskip+0.5cm

 \centerline{von Nikolay Moshchevitin\footnote
 {Steklow-Institut f\"ur Mathematik der Russischen Akademie der Wissenschaften. Diese Arbeit wurde unterst\"utzt durch RNF Grant No. 14-11-00433.
 }}
 
 \vskip+1cm

 \end{Large}
 
{\bf 1. Einleitung.}

Sei  $\alpha$
eine reelle Zahl und sei
$$
\psi_\alpha (t) = \min_{1\le q\le t, \, q\in \mathbb{Z}}
||q\alpha ||,\,\,\,\,\,\, t \ge 1
$$
die Funktion des Irrationalit\"atma\ss es f\"ur $\alpha$
(hier bezeichnet $||\xi || = \min_{a\in \mathbb{Z}}|\xi - a|
$ den Abstand zwischen $\xi$ und der n\"achstgelegenen ganzen Zahl). 
Viele diophantische Eigenschaften der Zahl  $\alpha$
k\"onnen 
durch Eigenschaften  der Funktion $
\psi_\alpha (t)$
ausgedr\"uckt werden.
Ziel dieser Arbeit ist es, 
einige Verallgemeinerungen  der Funktion   
$
\psi_\alpha (t)$
zu betrachten 
und einige entsprechende diophantische Behauptungen zu beweisen.

 {\bf 1.1. Kettenbr\"uche  und die Funktionen des Irrationalit\"atsma\ss es.}
 
 F\"ur eine  reelle irrationale Zahl  $\alpha \in [0,1]$
 definieren wir die unendliche Sequenz der Teilnenner
  $$\hbox{\got A} : { a}_1, { a}_2, a_3,...,{ a}_{n},...\,\,\, ,$$
  sodass
  $$
  \alpha =
  [0;a_1,a_2,a_3\dots,a_n,...]= \frac{1}{\displaystyle{a_1+\frac{1}{\displaystyle{a_2 +
\frac{1}{\displaystyle{a_3 +\dots +
\displaystyle{\frac{1}{a_n +\dots}} }}}}}} 
 .
 $$
 Wir betrachten die Naherungsbr\"uche
 $$
 \frac{p_n}{q_n} = [0;a_1,a_2,a_3\dots,a_n], \,\,\,\,\, n =1,2,3,...
 $$
 und setzen
 $$ \xi_n = ||q_n\alpha|| =|q_n\alpha - p_n| $$
und
 $$
 p_0=0, q_0 =1, \,\,\, p_{-1}= 1, q_{-1}=0.
 $$
 
 F\"ur die Gitterpunkte 
  $
  \pmb{z}_n = (p_n,q_n)\in\mathbb{Z}^2$
 gilt 
  \begin{equation}\label{1}
  \pmb{z}_{n} = a_n \pmb{z}_{n-1} + \pmb{z}_{n-2}.
 \end{equation}
 Seien
 $$
 \alpha_n = [a_n;a_{n+1},a_{n+2}, a_{n+3},...],\,\,\,\,\,\,\,
 \alpha_n^* = [0; a_n,a_{n-1},...,a_{1}],
 $$
 dann gelten
 die Gleichungen
   \begin{equation}\label{2p}
 \frac{q_{n-1}}{q_n} =\alpha_{n}^* 
 \end{equation}
und 
  \begin{equation}\label{2s}
 \frac{\xi_{n-1}}{\xi_n} =\alpha_{n+1} .
 \end{equation}

 Aus dem Lagrangeschen Gesetz der besten N\"aherung
folgt, dass
$$
\psi_\alpha (t) =\xi_n \,\,\,\text{f\"ur}\,\,\, q_n\le t < q_{n+1}.
$$
 Also die Funktion $\psi_\alpha (t) $ ist
 auf den Intervallen $(q_n, q_{n+1})$
 konstant. An den Stellen  $q_n$ 
 die Funktion $\psi_\alpha (t) $ 
 unstetig ist.

In dieser Arbeit  untersuchen wir  zwei Funktionen 
 $$
 \psi_\alpha^{[2]} (t)
 =\min_{
 \begin{array}{c}
 (q,p): \, q,p\in \mathbb{Z}, 1\le  q\le t, \cr
 (p,q) \neq (p_n, q_n) \,\forall\, n =0,1,2,3,...
 \end{array}
 }
 |q\alpha -p|
  $$
 und
 $$ \psi_\alpha^{[2]*} (t)
=\min_{
 \begin{array}{c}
 (q,p): \, q,p\in \mathbb{Z}, 1\le  q\le t, \cr
 p/q \neq p_n/q_n \,\forall\, n =0,1,2,3,...
 \end{array}
 }
 |q\alpha -p|,
 $$
 die
 mit  den ``zweitbesten Approximationen'' verbunden. 
  Beide Funktionen 
 $
 \psi_\alpha^{[2]} (t)
$
und 
 $
 \psi_\alpha^{[2]*} (t)
$
sind f\"ur
$
t\ge
1
$ 
definiert und sind stückweise konstant.  F\"ur diese Funktionen 
 betrachten wir die ganzzahligen Sequentzen der Stellen
 $$
 \hbox{\got Q} : \hbox{\got q}_1=1, \hbox{\got q}_2,...,\hbox{\got q}_{\nu},... 
 $$
 und
 $$
 \hbox{\got X}
 :
  \hbox{\got x}_1=1,
   \hbox{\got x}_2,....,  \hbox{\got x}_n,...
   $$ 
  wo diese Funktionen 
 unstetig sind.

Definieren wir den Wert
 $$
\hbox{\got t}_:=
\begin{cases} 
2,\,\,\text{falls entweder} \,\, a_1 =1, a_2\ge 2 \,\,
\text{oder}\,\,
 a_1 \ge 3
\cr
3,\,\,\text{falls entweder}\,\, a_1 = a_2 =1, a_3 \ge 2
\,\,
\text{oder}\,\,
a_1 =2, a_2\ge 2
\cr
4,\,\,\text{falls entweder}\,\, a_1=a_2=a_3 = 1
\,\,
\text{oder}\,\,
a_1 =2, a_2= 1
\end{cases}.
$$
 Man hat
 $$
 \psi_\alpha^{[2]} (t)
=
  \min_{
\begin{array}{c}
 1\le q\le t, 
 \, q\in \mathbb{Z}
 \cr q\neq q_n \,\forall\, n = 0,1,2,3,...
 \end{array}
 }
 ||q\alpha ||
 $$
 wenn  $t\ge \hbox{\got t}$.
  Also falls $\hbox{\got q}_\nu,
  \hbox{\got x}_\nu
  \ge \hbox{\got t}$,
  $$
 \text{
 f\"ur}  \,\,\,\,\hbox{\got q}_\nu \le t < \hbox{\got q}_{\nu+1}\,\,\,\, \text{ gilt } \,\,\,\,  \psi_\alpha^{[2]} (t) = || \hbox{\got q}_\nu\alpha ||,
 $$
 und
 $$
 \text{
 f\"ur} \,\,\,\, \hbox{\got x}_\nu \le t < \hbox{\got x}_{\nu+1} \,\,\,\, \text{ gilt }   \,\,\,\, \psi_\alpha^{[2]*} (t) = || \hbox{\got x}_\nu\alpha ||.
 $$ 
 Es ist klar,
 dass
 $$
 \psi_\alpha (t)
  <
  \psi_\alpha^{[2]} (t)
 \le \psi_\alpha^{[2]*} (t),\,\,\,\,\,
 \forall\, t.
 $$
  In den Punkten 2 und 3  ziehen wir die
  Regeln (S\"atze 
   1 und 2), wie
   die Sequenzen $ \hbox{\got Q},  \hbox{\got X}$   
   aus der Sequenz
   $ \hbox{\got A}$
   konstruieren   k\"onnen werden.

  {\bf 1.2. Das Lagrangesche Spektrum
 und die Spektra f\"ur die Funktionen  $\psi_\alpha^{[2]} $,   $\psi_\alpha^{[2]*}$.
 }

 F\"ur  irrationalen $\alpha$ bei
 $$
 \lambda (\alpha) = \liminf_{t\to \infty}\, t \psi_\alpha (t)
 $$
  bezeichnen
wir die Lagrangesche Konstant f\"ur $\alpha$.
 Die Menge 
 $$
 \mathbb{L}
=\{\lambda:\,\,\,\exists \alpha \in \mathbb{R}\setminus\mathbb{Q}\,\,\text{such that}\,\,
\lambda = \liminf_{t\to \infty} t \cdot\psi_\alpha (t)\}
 $$
 der Werte $\lambda (\alpha )$ hei\ss t
 das  Lagrangesche  Spectrum.

 Wir definieren  die Werte
 \begin{equation}\label{melichi}
 \hbox{\got j} (\alpha )=
\inf_{t\ge \hbox{\got t}} \, t
\,\psi_\alpha^{[2]} (t)
=
\inf_{\nu :\,\,\hbox{\got q}_\nu \ge \hbox{\got t}}\,  \hbox{\got k}_\nu
,\,
\,\,\,\,\,\,\,\,
\hbox{\got k}(\alpha )=
\liminf_{t\to \infty}\, 
\,t\psi_\alpha^{[2]} (t)
=
\liminf_{t\to \infty}\,\hbox{\got k}_\nu,
\end{equation}
wo
$$
\hbox{\got k}_\nu = \hbox{\got k}_\nu (\alpha) = 
\hbox{\got q}_\nu
|| \hbox{\got q}_\nu\alpha ||,
$$
und

    \begin{equation}\label{melichi}
 \hbox{\got j}^* (\alpha )=
\inf_{t\ge \hbox{\got t}}\, t
\,\psi_\alpha^{[2]*} (t)
=
\inf_{\nu :\,\,\hbox{\got x}_\nu \ge \hbox{\got t}}\,  \hbox{\got k}_\nu^*,
\,
\,\,\,\,\,\,\,\,
\hbox{\got k}^*(\alpha )=
\liminf_{t\to \infty}\, 
\,t\psi_\alpha^{[2]*} (t)
=
\liminf_{\nu\to \infty}\,\hbox{\got k}_\nu^*,
\end{equation}
wo
$$
\hbox{\got k}_\nu^* = \hbox{\got k}_\nu^* (\alpha) = 
\hbox{\got x}_\nu
|| \hbox{\got x}_\nu\alpha ||.
$$
 Es ist klar, dass
 $$
\hbox{\got j}(\alpha ) \le \hbox{\got k} (\alpha )
,
\,\,\,
\hbox{\got j}^*(\alpha ) \le \hbox{\got k}^* (\alpha )
.$$
Bemerken wir, dass f\"ur alle $\alpha$
 gilt
 \begin{equation}\label{impo1}
   \hbox{\got k}^* (\alpha )\ge  \hbox{\got k} (\alpha )
 \ge 2\lambda(\alpha ),
 \end{equation}
 und f\"ur
$\alpha\sim \sqrt{2}$ man hat 
\begin{equation}\label{impo2}
   \hbox{\got k}^* (\alpha ) =  \hbox{\got k} (\alpha )
 =2\lambda(\alpha ).
 \end{equation}
 Wir werden Formeln (\ref{impo1},\ref{impo2}) im Punkt 4.2 beweisen.

Wir untersuchen hier  die neuen Spektra
 $$
 \mathbb{L}_2
=\{\lambda:\,\,\,\exists \alpha \in \mathbb{R}\setminus\mathbb{Q}\,\,\text{mit}\,\,
\lambda = \liminf_{t\to \infty} t \cdot\psi_\alpha^{[2]} (t)\}
 $$
 und
 $$
 \mathbb{L}_2^*
=\{\lambda:\,\,\,\exists \alpha \in \mathbb{R}\setminus\mathbb{Q}\,\,\text{mit}\,\,
\lambda = \liminf_{t\to \infty} t \cdot\psi_\alpha^{[2]*} (t)\},
 $$
die
 durch die Funktionen $\psi^{[2]}(t)$ und  $\psi^{[2]*}(t)$
      definieren werden. 
   Wir formulieren and beweisen zwei S\"atze \"uber die Struktur der Spektra
  $\mathbb{L}_2$ und $\mathbb{L}_2^*$   in den Punkten 
  4.3 und 4.3.

 {\bf 1.3. Der Legendresche Satz und Verallgemeinerungen.}

 Die bekannteste Version des Satzes von Legendre ist die folgende:
 {\it
 Falls $\left|\alpha -\frac{p}{q}\right|<\frac{1}{2q^2}$ und
 $(p,q) =1$,
 der Bruch $\frac{p}{q}$ ein N\"acherungsbruch f\"ur $\alpha$ ist.}
  Aber das originale Ergebnis  von Legendre (siehe \cite{lege} und \cite{Pe} S.42-45, \cite{Ve} S. 42-43.)
   ist eine stärkere Behauptung. 
  Es gibt uns ein
Kriterium f\"ur $\alpha$, den Bruch $\frac{p}{q}$ als einen N\"acherungsbruch zu haben.
 
{\bf Satz von Legendre.}
\,\, {\it Sei $\frac{p}{q}, \, (p,q) =1$ eine rationale Zahl.
F\"ur irrationale $\alpha$ betrachten wir den Wert
$$
\theta = q \cdot (q\alpha - p).
$$
F\"ur $\frac{p}{q}$ definieren wir die Kettenbruch
$$
\frac{p}{q} = [a_0; a_1,...,a_t], 
$$
mit $a_t \ge 1$,
wo 
$t$ gerade ist falls $\theta <0$,
und $t$ ungerade ist falls $\theta >0$. Definieren wir den Bruch
 $$
\frac{p'}{q'} = [a_0; a_1,...,a_{t-1}]. 
$$
Dann die notwendige und hinreichende Bedingung daf\"ur, dass
$
\frac{p}{q} $
 ein N\"acherungsbruch von $\alpha$ ist, lautet daher:
 $$
 |\theta| < \frac{q}{q+q'}
 .
 $$}
 
 Lucas (siehe \cite{LU} S. 447-449 und \cite{lehmer} S. 229)\footnote{Der Autor  dankt Prof. Yu. Nesterenko f\"ur diese Referenzen.} hat ein gleiches 
 Kriterium formuliert, mit keine Referenz. Wir formulieren die Version  des Kriteriums aus \cite{lehmer}.

 {\bf Satz von Lucas.}
\,\, {\it Seien $\frac{p_{t-1}}{q_{t-1}}$
und 
$\frac{p_{t}}{q_{t}}$
zwei sukzessive N\"acherungsbr\"uche f\"ur eine Zahl $\eta$. Dann 
$\frac{p_{t-1}}{q_{t-1}}$
und 
$\frac{p_{t}}{q_{t}}$ sind
zwei sukzessive N\"acherungsbr\"uche f\"ur $\alpha$ dann und nur dann, wenn
$$
\left|\alpha - \frac{p_t}{q_t}\right| <\frac{1}{q_t(q_t+q_{t-1})}
.
$$}

Die Natur dieser S\"atze ist klar. Gegeben einen Bruch $\frac{p}{q}$,  sollen wir die Farey-Folge $q$-ter Ordnung
$\hbox{\got F}_q$ betrachten. Seien  $r_{j-1}< r_j = \frac{p}{q}< r_{j+1}$ drei sukzessive Elemente aus $\hbox{\got F}_q$.
F\"ur den einzigen Kettenbruch
$$
\frac{p}{q} = [a_0;a_1,..., a_t]
,\,\, a_t\ge 2
$$
 haben wir 
$$
r_{j-1}= \frac{p_-}{q_-} = [a_0;a_1,...,a_{t-1}] ,\,\,\,
r_{j+1} =  \frac{p_+}{q_+}=[a_0;a_1,...,a_{t}-1] ,
$$
{oder}
$$
r_{j-1}= \frac{p_-}{q_-} = [a_0;a_1,...,a_{t}-1] ,\,\,\,
r_{j+1} =  \frac{p_+}{q_+}=[a_0;a_1,...,a_{t-1}] ,
$$
bez\"uglich der Parit\"at. Dann
$
\frac{p}{q} $
 ein N\"acherungsbruch von $\alpha$ ist dann und nur dann, wenn
 $$
 \alpha \in \left(
 \frac{p+p_-}{q+q_-}, 
 \frac{p+p_+}{q+q_+}
 \right).
  $$
  
  Bemerken wir, dass das Kriterium von Legendge ist eine Aussage, die gibt  uns f\"ur  
   eine gegebene rationale Zahl $\frac{p}{p}$ eine notwenige und hinreichende  Bedingung f\"ur
   $\alpha$, den Bruch $\frac{p}{q}$ als einen  N\"acherungsbruch zu haben. Wir wollen hier zwei Behauptungen formulieren, die 
   f\"ur eine gegebene irrationale Zahl $\alpha$  geben die Bedingungen f\"ur den Bruch $\frac{p}{q}$ ein  N\"acherungsbruch f\"ur $\alpha$ zu sein.
 Diese Behauptungen  folgen sofort aus den
 Definitionen 
 der Funktionen   $\psi_\alpha^{[2]} (t)$ und   $\psi_\alpha^{[2]*} (t)$:

  \vskip+0.4cm
 \noindent 
 { 1)}
 {\it
  Falls 
  $$
  ||q\alpha|| <
  \psi_\alpha^{[2]} (q),
  $$ 
  $q$ ein Nenner
  eines N\"acherungsbruches f\"ur $\alpha$ ist.}

   \vskip+0.4cm
  
   \noindent 
 { 2)}
 {\it
  Falls 
  $$
  |q\alpha -p| <
  \psi_\alpha^{[2]*} (q),
  $$ 
   und
 $(p,q) =1$,
 der Bruch $\frac{p}{q}$ ein N\"acherungsbruch f\"ur $\alpha$ ist.}
    \vskip+0.4cm

 Aus den Definitionen
 der Werte
 $
  \hbox{\got j} (\alpha),
   \hbox{\got k} (\alpha)
 $
 und
  $
  \hbox{\got j}^* (\alpha),
   \hbox{\got k}^* (\alpha)
 $
 haben wir die folgenden Behauptungen:

  \vskip+0.4cm
 \noindent 
 {3)}
 {\it
   Falls  $ q ||q\alpha|| <  \hbox{\got j}(\alpha)$,
  $q$ ein Nenner
  eines N\"acherungsbruches f\"ur $\alpha$ ist.}

   \vskip+0.4cm
  
   \noindent 
 {4)}
 {\it
 F\"ur $\varepsilon >0$  gibt es eine effektiv berechenbare Konstante
 $ T = T(\alpha, \varepsilon)$ mit  der folgenden Eigenschaft:
 falls $ q\ge T $ und $ q 
    ||q\alpha ||\le
 \hbox{\got k}(\alpha)
 -\varepsilon $,
$q$ ein Nenner
  eines N\"acherungsbruches f\"ur $\alpha$ ist.}

   \vskip+0.4cm
   
   \noindent 
  {5)}
  {\it
 Falls 
  $ q\cdot
  |q\alpha -p| <
 \hbox{\got j}^*(\alpha)$ 
   und
 $(p,q) =1$,
 der Bruch $\frac{p}{q}$ ein N\"acherungsbruch f\"ur $\alpha$ ist.}

   \vskip+0.4cm
   
   \noindent 
   {6)}
    {\it
 F\"ur $\varepsilon >0$  gibt es eine effektiv berechenbare Konstante
 $ T^* = T^*(\alpha, \varepsilon)$ mit  der folgenden Eigenschaft:
 falls $ q\ge T^* $ und $ q \cdot
    |q\alpha -p|\le
 \hbox{\got k}(\alpha)-\varepsilon,\, (q,p) = 1
  $,
 der Bruch $\frac{p}{q}$ ein N\"acherungsbruch f\"ur $\alpha$ ist.}

   \vskip+0.4cm

 {\bf 1.4. Die Beispiele.}

  Diskutieren wir einige Beispiele f\"ur speziele Zahlen.
  
   \begin{itemize}

 \item
 Betrachten wir die Zahl
 $\tau =\frac{1+\sqrt{5}}{2}$. Dann
 $$
 \hbox{\got  j} (\tau) =   \hbox{\got  j}^* (\tau) = 8(\sqrt{5}-2),\,\,\,\,\,\,
  \hbox{\got  j} (\tau) = \frac{4}{\sqrt{5}},\,\,\,\,\,\,
    \hbox{\got  k}^* (\tau) = \sqrt{5}.
    $$
    Also
 
 \noindent
 1) falls
 $q||q\tau || < 8(\sqrt{5}-2)$,
 $q$ eine Fibonacci-Zahl ist;
 
 \noindent
 2) zu jedem $\varepsilon>0$, falls 
 $q$ hinreichend gro\ss \, ist und
 $||q\tau||\le \frac{4 -\varepsilon}{\sqrt{5}\cdot q}$,
 die Zahl
$q$ eine Fibonacci-Zahl ist;

\noindent
3)
zu jedem $\varepsilon>0$ falls 
 $q$ hinreichend gro\ss \, ist und
  $\left|\tau -\frac{p}{q}\right|\le \frac{\sqrt{5}-\varepsilon}{\cdot q^2}, (p,q) = 1$,
  der Bruch
 $\frac{p}{q}$ ein 
 N\"acherungsbruch
 f\"ur
 $\tau$ ist;
 
 \item
  Betrachten wir die Zahl 
 $\xi =\frac{1+\sqrt{17}}{2}$. Dann $\hbox{\got k} (\xi ) = \frac{4}{\sqrt{17}}$
 und
 zu jedem $\varepsilon >0$ falls  $q$
 hinreichend gro\ss \, ist und
  $||q\tau||\le \frac{4-\varepsilon}{\sqrt{17}\cdot q}$,
  die Zahl $q$
 ein
  Nenner
  eines N\"acherungsbruches f\"ur $\xi$ ist.

\item
  Betrachten wir die Zahl 
 $e =\sum_{k=0}^\infty \frac{1}{k!}$.
 Dann $\hbox{\got k}^* (e ) = \frac{3}{2}$
 und
 zu jedem $\varepsilon >0$ falls  $q$
 hinreichend gro\ss \, ist und
   $\left|e -\frac{p}{q}\right|\le \left(\frac{3}{2}-\varepsilon\right)\cdot \frac{1}{\cdot q^2}$ mit $(p,q) = 1$,
 der Bruch 
 $\frac{p}{q}$
 ein
 N\"acherungsbruch
 f\"ur
 $e$ ist.
\end{itemize}

{\bf 2.  \"Uber die Funktion  $\psi_\alpha^{[2]} (t)$.}

Wir formulieren hier eine allgemeine
Regel, die die Sequenz $ \hbox{\got Q}$ aus der Sequenz $\hbox{\got A}$ konstruirt.

 {\bf Satz 1.}
 \,\,
 {\it Die Sequenz $\hbox{\got Q}$ wird mit der folgenden Regel  aus der Sequenz $\hbox{\got A}
 $ erhalten:

\noindent
 1. jedes Element $ a_n \ge 3$ wird durch 
  sukzessive Zahlen
 \begin{equation}\label{ai1}
q_{n-2}+q_{n-1},\,\,\, 
2q_{n-1},\,\,\,
q_{n}-q_{n-1}
\end{equation}
ersetzen;
 
\noindent
 2. jedes Element $ a_n =2$ wird  durch 
 eine Zahl 
 $$
 q_{n}-q_{n-1}
 $$
 ersetzen;
 
 \noindent
 3. falls $a_{n-1}\neq 1, a_{n}=1, a_{n+1}\ge 2$,
 das Element $a_n =1$
 wird durch 
 eine Zahl $$q_{n-2}+q_n = 2q_{n-2}+q_{n-1}$$ ersetzen;
 
 \noindent
4. falls $r\ge 2$ und $a_{n-1}\neq 1, a_{n}= ...=a_{n+r-1}=1, a_{n+r}\ge 2$,
die Elemente $ a_{n}= ...=a_{n+r-1}=1$ werden
durch 
 sukzessive Zahlen
\begin{equation}\label{ai2}
 2q_{n-2}+q_{n-1}, 2q_{n-1}, 2q_n,...,2q_{n+r-2},  2q_{n+r-3}+q_{n+r-2}
 \end{equation}
ersetzen;

 \noindent
5. falls  $a_{n-1}\neq 1, a_{j}=1,\,\forall \, j \ge n  $,
die Elemente $ a_{n} = a_{n+1} = a_{n+2}=...$ werden
durch 
 Zahlen
$$
 2q_{n-2}+q_{n-1}, 2q_{n-1}, 2q_n, 2q_{n+1},...
$$
ersetzen.

}
 
 {\bf Bemerkung 1.}
 
 {\it
  \noindent
  1.1.
 Falls $n\ge 3$ oder $ n =2, a_n \ge 3$, alle Zahlen in den Formeln (\ref{ai1}) und (\ref{ai2})
  unterschiedlich sind.
 So werden in (\ref{ai1}) drei verschiendene Zahlen geschrieben, und in (\ref{ai2}) werden $r+2$
 verschiendene Zahlen geschrieben.

  \noindent
  1.2.
 F\"ur $n = 1$ und $a_1=3$ in
 (\ref{ai1}) haben wir  nur zwei  verschiedene Zahlen
 $$
1= q_{-1}+q_0 , \,\,\, 2 = 2q_0 = q_1-q_0.
$$

 \noindent
1.3.
 F\"ur $n = 1$, $a_1=1$ und $ r \ge 4$ in
 (\ref{ai2}) haben wir  nur $r+1$ verschiedene Zahlen
$$
1= 2q_{-1}+q_0,
2= 2q_0 = 2q_1, 4 = 2q_2, ... ,  
2q_{r-1},  2q_{r-2}+q_{r-1}.
$$

 \noindent
1.4.
 F\"ur $n = 1$, $a_1=1$ und $ r =3$ in
 (\ref{ai2}) haben wir  drei verschiedene Zahlen
$$
1= 2q_{-1}+q_0,
2= 2q_0 = 2q_1, 4 =    
2q_{2}= 2q_{1}+q_{2}.
$$

 \noindent
1.5.
 F\"ur $n = 1$, $a_1=1$ und $ r =2$ in
 (\ref{ai2}) haben wir  drei verschiedene Zahlen
$$
1= 2q_{-1}+q_0,
2= 2q_0 = 2q_1, 
3=  2q_{0}+q_{1}
$$

\noindent
2. Falls $ a_1=a_2=1, a_3 \ge 2$, 
f\"ur  die Teilnenner $ a_1=a_2=1$
Punkt 4 des Satzes 1 gibt   das Element 
$3 = 2q_0+q_1$ der Sequenz $\hbox{\got Q}$. F\"ur den  Teilnenner $a_3\ge 3$  Punkt 1 des Satzes 1 
(oder Punkt 2, falls $ a_3 =2$)
gibt das Element 
$ 3= q_1+q_2$ auch.
In allen anderen F\"allen
f\"ur  verschiedene $n$ gibt
die Regel des Satzes 1  verschiendene Zahlen in $\hbox{\got Q}$.
 }

F\"ur den Gitterpunkt $\pmb{z} =(q,p) \in \mathbb{Z}^2$  
setzen
wir
$$
q = q(\pmb{z}),\,\,\, \xi =\xi (\pmb{z} ) = |q\xi - p|. 
$$

Nun formulieren wir eine einfache und wichtige Behauptung, die aus der
Definition der Funktion $\psi^{[2]}(t)$ folgt sofort.

{\bf Haupthilfssatz 1.}\,\,{\it
Um zu zeigen, dass 
$x_1< x_2$  zwei sukzessive Elemente aus $\hbox{\got Q}$ sind,
genügt es, folgendes zu beweisen:

\noindent
$\bullet$ f\"ur einige  $y_1,y_2 \in \mathbb{Z}$ die Gitterpunkte
$\pmb{w}_1 = (x_1,y_1), 
\pmb{w}_2 = (x_2,y_2)
$
auf dem Rand des Parallelogrammes
$$ 
\Pi [\pmb{w}_1, 
\pmb{w}_2] =
\{
\pmb{z} = (x,z):\,\,
0\le x\le x_2,\,\,
|\alpha x - y| \le \xi (\pmb{w}_1) = |\alpha x_1-y_1|\}
$$
liegen;

\noindent
$\bullet$
das Parallelorgamm 
$$
\Omega [\pmb{w}_1]  =
\{
\pmb{z} = (x,z):\,\,
0\le x\le x_1,\,\,
|\alpha x - y| \le \xi (\pmb{w}_1) = |\alpha x_1-y_1|\}
\subset \Pi [\pmb{w}_1, 
\pmb{w}_2] 
$$
hat einen inneren Punkt $\pmb{z}_{\nu_1}  = (p_{\nu_1}, q_{\nu_1})$;
 
 \noindent
$\bullet$
das Parallelorgamm 
$$
\Omega [
\pmb{w}_2] =
\{
\pmb{z} = (x,z):\,\,
0\le x\le x_2,\,\,
|\alpha x - y| \le \xi (\pmb{w}_2) = |\alpha x_2-y_2|\}
\subset \Pi [\pmb{w}_1, 
\pmb{w}_2] 
$$
hat einen inneren Punkt $\pmb{z}_{\nu_2} = (p_{\nu_2}, q_{\nu_2})$;
 
\noindent
$\bullet$
falls $\pmb{z} \in \Pi [\pmb{w}_1, 
\pmb{w}_2] $
ein Gitterpunt ist, dann
$\pmb{z} \in \{\pmb{0}, \pmb{w}_1,\pmb{w}_2\}$ oder
$\pmb{z} =\pmb{z}_\nu$ mit einem $\nu\ge -1$.}

In das Folgendes,
 Elemente der Sequentz $\hbox{\got Q}$ entsprechen roten Punkten auf den Figuren 1 - 7.
Wir leiten nun aus diesem Haupthilfssatz Hilfss\"atze 1 - 7  her.
Diese Hilfss\"atze  beweisen die Regeln aus Satz 1.
Sie
erhalten 
aus den Teilnenner f\"ur $\alpha$
die
Elemente der Sequenz $\hbox{\got Q}$  
und
 verbinden die einzelnen Bl\"ocke. 
 Satz 1 folgt daraus.

{\bf Hilfssatz 1.}
\,\,{\it
Sei $a_n \ge 3$. Dann 
$$
q_{n-2}+q_{n-1},\,\,\,
2q_{n-1},\,\,\,
q_{n}-q_{n-1}$$
sind drei sukzessive Elemente der Sequenz $\hbox{\got Q}$.
}

\begin{figure}[h]
  \centering
  \begin{tikzpicture}[scale=1.6]

  \draw[color=black] (0,0) -- (7.7,0);
         \draw (7.35,-0.2) node[]{\begin{small}$y=\alpha x$\end{small}};
  
  \draw[color=black] (0,-2.1) -- (0,3.1);

    \draw[color=black] (-0.4,-1.4) -- (0.8,2.8);
    \draw[color=black] (0.8,2.8) -- (8.6,1.37);
    \draw[color=black] (0.4,1.4) -- (8.2,-0.03);
 \draw[color=black] (0,0) -- (4.2,-0.77);
  \draw[color=black] (-0.4,-1.4) -- (2.6,-1.94);
 
    \draw[color=black] (0.8,-1.62) -- (2,2.6);
  
    \draw[color=black] (2,-1.84) -- (3.2,2.36);

   \node[fill=black,circle,inner sep=1.2pt]   at (4,0.74) {};
   \node[fill=black,circle,inner sep=1.2pt]   at (5.2,0.52) {};
  
   \node[fill=black,circle,inner sep=1.2pt]   at (7.6,0.08) {};
   
    \node[fill=black,circle,inner sep=1.2pt]   at (0,0) {};
     \node[fill=black,circle,inner sep=1.2pt]   at (-0.4,-1.4) {};
    \node[fill=black,circle,inner sep=1.2pt]   at (0.4,1.4) {};   
       \node[fill=black,circle,inner sep=1.2pt]   at (0.8,2.8) {};
         \node[fill=black,circle,inner sep=1.2pt]   at (0.8,-1.62) {};
           \node[fill=black,circle,inner sep=1.2pt]   at (2,-1.84) {};

    \node[fill=black,circle,inner sep=1.2pt]   at (2,2.58) {};
    \node[fill=black,circle,inner sep=1.2pt]   at (3.2,2.36) {};
       \node[fill=black,circle,inner sep=1.2pt]   at (4.4,2.14) {};
           \node[fill=black,circle,inner sep=1.2pt]   at (5.6,1.92) {};
               \node[fill=black,circle,inner sep=1.2pt]   at (6.8,1.7) {};
      \node[fill=black,circle,inner sep=1.2pt]   at (8,1.48) {};
      \node[fill=black,circle,inner sep=1.2pt]   at (1.2,-0.22) {};

      \node[fill=black,circle,inner sep=1.2pt]   at (3.6,-0.66) {};   
      \node[fill=black,circle,inner sep=1.2pt]   at (2.8,0.96) {};   
  
  \node[draw=black,circle,inner sep=2pt] at (0.4,1.4) {};
  \draw (0.75,1.5) node[]{$\pmb{z}_{n-2}$};

  \node[draw=black,circle,inner sep=2pt] at (7.6,0.08) {};
   \draw (7.8,0.2) node[]{$\pmb{z}_{n}$};

  \node[draw=black,circle,inner sep=2pt] at (1.2,-0.22) {};
   \draw (1.55,-0.15) node[]{$\pmb{z}_{n-1}$};
  
         \draw (2.8,-1.9) node[]{\begin{small}$\ell_{-1}$\end{small}}; 
  
  \draw (4.3,-0.7) node[]{\begin{small}$\ell_{0}$\end{small}};
  
 \draw (8.25,0.1) node[]{\begin{small}$\ell_{1}$\end{small}}; 
  
    \draw (8.65,1.5) node[]{\begin{small}$\ell_{2}$\end{small}};

    \draw[color=black] (0,1.18) -- (2.4,1.18);
    \draw[color=black] (0,-1.18) -- (1.7,-1.18);
    \draw[color=black] (2.1,-1.18) -- (2.4,-1.18);
    \draw[color=black] (2.4,1.18) -- (2.4,-1.18);
        \draw (1.9,-1.18) node[]{$\Pi_1$};

    \draw[dashed,color=black]  (0,1.16) -- (1.6,1.16);
    
      \draw (1.57,-0.76) node[]{${\Omega}_1$};

    \draw[dashed,color=black]  (0,-1.16) -- (1.6,-1.16);
    \draw[dashed,color=black]  (1.6,-1.16) -- (1.6,-0.89);
    \draw[dashed,color=black]  (1.6,-0.57) -- (1.6,1.16);

    \draw[color=black] (0,0.44) -- (6.4,0.44);
    \draw[color=black] (0,-0.44) -- (4.95,-0.44);
    \draw[color=black] (5.3,-0.44) -- (6.4,-0.44);
    \draw[color=black] (6.4,0.44) -- (6.4,-0.44);
        \draw (5.15,-0.44) node[]{$\Pi_2$};

    \draw[dashed,color=black]  (0,0.3) -- (4,0.3);
    \draw[dashed,color=black]  (4.3,0.3) -- (6.4,0.3);
      \draw (4.2,0.26) node[]{${\Omega}_3$};

    \draw[dashed,color=black]  (0,-0.3) -- (6.38,-0.3);
    \draw[dashed,color=black]  (6.38,-0.3) -- (6.38,0.3);

     \draw (2.35,1.35) node[]{$\pmb{z}_{n-2}+\pmb{z}_{n-1}$};  

          \draw (2.8,-0.65) node[]{$2\pmb{z}_{n-1}$};  
    
         \draw (6.97,0.44) node[]{$\pmb{z}_{n}-\pmb{z}_{n-1}$};  
    
  \node[draw=black,fill=red,circle,inner sep=2pt] at (1.6,1.18) {};
  \node[draw=black,fill=red,circle,inner sep=2pt] at (2.4,-0.44) {}; 
  \node[draw=black,fill=red,circle,inner sep=2pt] at (6.4,0.3) {}; 
  
  \end{tikzpicture}
  
   \caption{$a_n =6$} \label{fig1}
\end{figure}

Beweis (siehe Fig. 1).
\, Beweisen wir, dass
\begin{equation}\label{d}
 \xi (\pmb{z}_{n-2}+\pmb{z}_{n-1}) =\xi_{n-2}-\xi_{n-1}
 >
 \xi (2\pmb{z}_{n-1}) =
 2\xi_{n-1}
 >
 \xi (\pmb{z}_n - \pmb{z}_{n-1} ) =
 \xi_{n-1}+\xi_n.
 \end{equation}
Wegen $a_n \ge 3$  hat man
$ \frac{\xi_{n-2}}{\xi_{n-1}} =\alpha_n >a_3 \ge 3$, und die erste Ungleichung aus (\ref{d}) folgt.
Wegen $ \xi_{n-1}>\xi_n$ haben wir die zweite Ungleichung bewiesen auch.
 
  Betrachten wir die Parallelogramme
$$
\Pi_1 = \Pi[\pmb{z}_{n-2}+\pmb{z}_{n-1}, 2\pmb{z}_{n-1}]= 
\{
(x,y) \in \mathbb{R}^2:\,\, 
0\le x\le  2q_{n-1},\,\,\,
|\alpha x-y|\le \xi_{n-2}-\xi_{n-1}\},
$$
$$
\Pi_2 =
\Pi[ 2\pmb{z}_{n-1},\pmb{z}_n-\pmb{z}_{n-1}] = 
\{
(x,y) \in \mathbb{R}^2:\,\,
0\le x\le  q_n- q_{n-1},\,\,\,
|\alpha x-y|\le 2\xi_{n-1}\},
$$
und
$$
\Omega_1 =\Omega 
[
 {z}_{n-2}+\pmb{z}_{n-1}]
 \subset \Pi_1 ,
\,\,\,
 {\Omega}_2=
\Omega_2 [2\pmb{z}_{n-1}]
=
\Pi_1\cap\Pi,
\,\,\,
\Omega_3 =\Omega[\pmb{z}_n-\pmb{z}_{n-1}]
  \subset \Pi_2 .
$$

Jedes Parallelogramm $\Pi_j$  hat zwei Gitterpunkte auf dem Rand.
Die Punkte $ \pmb{z}_{n-2}+\pmb{z}_{n-1}$ und $2\pmb{z}_{n-1}$ sind
diese Punkte f\"ur $\Pi_1$.
Die Punkte $ 2\pmb{z}_{n-1}, \pmb{z}_n - \pmb{z}_{n-1}$ sind
diese Punkte f\"ur $\Pi_2$.
Der Punkt $2\pmb{z}_{n-1}$ geh\"ort  beiden Parallellogrammen $\Pi_1,\Pi_2$.

Es ist klar, dass
\begin{equation}\label{k1}
q_{n-1} <  q_{n-2}+ q_{n-1}<
\min \{ 2q_{n-1}, q_n\}
.
\end{equation}
Aus (\ref{1})
folgt
$$
\xi_n = \xi_{n-2}- a_n \xi_{n-1}
,$$
 sodass
 \begin{equation}\label{k2}
 \xi_{n-1}<
  \xi (\pmb{z}_{n}-\pmb{z}_{n-1}) =
 \xi_{n}+\xi_{n-1}
 <
  \xi (\pmb{z}_{n-2}+\pmb{z}_{n-1}) =
 \xi_{n-2}-\xi_{n-1}
 < 2\xi_{n-1}. 
 \end{equation}
Wir sehen aus (\ref{k1},\ref{k2}), dass der Gitterpunkt $\pmb{z}_{n-1}$
ein inner Punkt f\"ur alle Parellelogramme $\Pi_1,\Pi_2,{\Omega}_1,{\Omega}_2, \Omega_3$ ist:
\begin{equation}\label{i}
 \pmb{z}_{n-1} \in {\rm int}\, ({\Omega}_1\cap {\Omega}_3) \subset
  {\rm int}\,
 (\Omega_2).
\end{equation}
Nun beweisen wir, dass
\begin{equation}\label{2}
(\Pi_1\cup\Pi_2) \cap\mathbb{Z}^2 =\{\pmb{0},
 \pmb{z}_{n-2}+\pmb{z}_{n-1}, \pmb{z}_{n-1}, 2\pmb{z}_{n-1}, \pmb{z}_n - \pmb{z}_{n-1}\}.
 \end{equation}
 Das Gitter $\mathbb{Z}^2$ teilt sich in  den Linien
$$
\mathbb{Z}^2 =
\bigcup_{u \in \mathbb{Z}} \,\ell_u,\,\,\,\
\ell_u =\{
\pmb{z}\in \mathbb{Z}^2:
\,\,
\pmb{z} = u\pmb{z}_{n-2}+v\pmb{z}_{n-1}, \,\, v \in \mathbb{Z}\}.
$$

Falls $\pmb{z}\in \ell_u$ mit $ u \neq 0,1$,  
$$
\xi (\pmb{z})>\xi_{n-2}>  \xi_{n-2}-\xi_{n-1}
=
\xi (\pmb{z}_{n-2}+\pmb{z}_{n-1})
,
$$
und $ \pmb{z} \not\in \Pi_1\cup\Pi_2 $.

Falls $\pmb{z}\in \ell_0$, haben wir
$\pmb{z}_{n-1}, 2\pmb{z}_{n-1}\in  \Pi_1\cap\Pi_2$, aber 
f\"ur $ v \ge 3$ gilt
$$
\xi (v\pmb{z}_{n-1}) \ge  3 \xi_{n-1} 
 > \xi_{n}+\xi_{n-1}=
 \xi (\pmb{z}_{n}-\pmb{z}_{n-1}),\,\,\,
 q(v\pmb{z}_{n-1}) > 2q_{n-1},
 $$
 und $ \pmb{z} \not\in \Pi_1\cup\Pi_2 $.
 
 Falls $\pmb{z}
 =
 \pmb{z}_{n-2}+v \pmb{z}_{n-1}
 \in \ell_1$,  haben wir  
$\pmb{z}_{n-2}+\pmb{z}_{n-1} \in  \Pi_1$ und $ \pmb{z}_n -\pmb{z}_{n-1} \in\Pi_2$, aber
f\"ur $ 2\le v \le a_{n} -2$ gilt
$$
 \xi (v\pmb{z})\ge
 2\xi_{n-1}+\xi_n > 2\xi_{n-1} = \xi (2\pmb{z}_{n-1})
 ,\,\,
  q(\pmb{z}) = q_{n-2} + vq_{n-1} > 2q_{n-1},
 $$
 und $ \pmb{z} \not\in \Pi_1\cup\Pi_2 $. F\"ur $ v \ge a_n$ alles ist klar,
 und (\ref{2}) bewiesen ist.
 
  Hilfssatz 1 folgt aus (\ref{d},\ref{2}) und (\ref{i}). 
 $\Box$

 {\bf Hilfssatz 2.}
\,\,{\it
Sei $a_n \ge 2$ und $a_{n+1}\ge 2$. Dann 
$$
q_{n}-q_{n-1},\,\,\,
q_{n}+q_{n-1}$$
sind zwei sukzessive Elemente der Sequenz $\hbox{\got Q}$.
}

\begin{figure}[h]
  \centering
  \begin{tikzpicture}[scale=1.3]
    
    \draw[color=black] (-1.2,-1.6) -- (2.4,3.2);
    \draw[color=black] (0,0) -- (2.6,-0.4);
    \draw[color=black] (1.4,-2) -- (3.8,1.2);
 \draw[color=black] (-1.2,-1.6) -- (1.4,-2);
  \draw[color=black] (1.2,1.6) -- (3.8,1.2);

 \node[fill=black,circle,inner sep=1.2pt]   at (0,0) {};
 \node[fill=black,circle,inner sep=1.2pt] at (-1.2,-1.6) {};
  \node[fill=black,circle,inner sep=1.2pt]   at (1.2,1.6) {};
  \node[fill=black,circle,inner sep=1.2pt]   at (2.4,3.2) {};
  \node[fill=black,circle,inner sep=1.2pt] at (2.6,-0.4) {};

  \node[draw=black,circle,inner sep=2pt] at (1.2,1.6) {};
  \node[draw=black,circle,inner sep=2pt] at (2.6,-0.4) {};

  \draw[color=black] (0,0) -- (5,0);
         \draw (4.6,0.1) node[]{\begin{small}$y=\alpha x$\end{small}};
  
  \draw[color=black] (0,-2.3) -- (0,3);

    \draw[dashed,color=black]  (0,1.2) -- (3.8,1.2);

    \draw[dashed,color=black]  (0,-1.2) -- (2.85,-1.2);
    \draw[dashed,color=black]  (3.15,-1.2) -- (3.78,-1.2);
    
    \draw[dashed,color=black]  (3.77,1.2) -- (3.77,-1.2);
 
    \draw[color=black]  (3.8,2) -- (3.8,-2); 

    \draw[color=black] (0,2) -- (2.8,2);
   \draw[color=black] (3.15,2) -- (3.8,2);
  
  \draw[color=black] (0,-2) -- (3.8,-2);
 
 \draw[dashed,color=black]  (1.4,1.97) -- (1.4,0.77);
  \draw[dashed,color=black]  (1.4,0.5) -- (1.4,-2);

 \draw[dashed,color=black]  (0,1.97) -- (1.4,1.97);
 \draw[dashed,color=black]  (0,-1.97) -- (1.4,-1.97);
 
  \draw (0.95,1.75) node[]{\begin{small}$\pmb{z}_{n-1}$\end{small}};
  \draw (4.5,1.3)node[]{\begin{small}$\pmb{z}_{n-1}+\pmb{z}_n$\end{small}};
 \draw (2.9,-0.4)node[]{\begin{small}$\pmb{z}_{n}$\end{small}};
 \draw (2.1,-2.15)node[]{\begin{small}$\pmb{z}_{n}-\pmb{z}_{n-1}$\end{small}};

  \draw (3,2.03)node[]{$\Pi$};
  \draw (1.46,0.62)node[]{$\Omega_1$}; 
  \draw (3,-1.2)node[]{$\Omega_2$};
 
  \node[draw=black,fill=red,circle,inner sep=2pt] at (1.4,-2) {};
  \node[draw=black,fill=red,circle,inner sep=2pt] at (3.8,1.2) {}; 
  \end{tikzpicture}
  
   \caption{$a_n, a_{n+1} \ge 2$} \label{fig2}
\end{figure}

Beweis  (siehe Fig. 2).\,
 Bemerken wir, dass
 die Br\"uche
 $ \frac{p_{n}-p_{n-1}}{q_{n}-q_{n-1}},
  \frac{p_{n}-p_{n-1}}{q_{n}-q_{n-1}}$ 
  sind N\"aherungsbr\"uche nicht, und
 $$
 \xi (\pmb{z}_n-\pmb{z}_{n-1})=
 \xi_n+\xi_{n-1}> \xi_n -\xi_{n-1} =
  \xi (\pmb{z}_n+\pmb{z}_{n-1}).
  $$
  Betrachen wir das Parallelogramm  
 $$
 \Pi = \Pi [\pmb{z}_n-\pmb{z}_{n-1},\pmb{z}_n+\pmb{z}_{n-1}]
=\{
(x,y) \in \mathbb{R}^2:\,\,
0\le x\le  q_{n}+q_{n-1},\,\,\,
|\alpha x-y|\le \xi_{n-1}+\xi_{n}\},
$$
und die Parallelogramme
$$
\Omega_1 = \Omega[\pmb{z}_n - \pmb{z}_{n-1}]=
\{
(x,y) \in \mathbb{R}^2:\,\,
0\le x\le  q_{n}-q_{n-1},\,\,\,
|\alpha x-y|\le \xi_{n-1}+\xi_{n}\}\subset \Pi,
$$
$$
\Omega_2 = 
\Omega[\pmb{z}_n+  \pmb{z}_{n-1}]
=\{
(x,y) \in \mathbb{R}^2:\,\,
0\le x\le  q_{n}+q_{n-1},\,\,\,
|\alpha x-y|\le \xi_{n-1}-\xi_{n}\}\subset\Pi.
$$
Die Punkte $\pmb{z}_n \pm \pmb{z}_{n-1} $
liegen auf dem Rand des Parallelogrammes $\Pi$.

Wegen
$ a_n \ge 2$
haben wir
\begin{equation}\label{l1}
q_{n-1}< q_n -q_{n-1}.
\end{equation}
Es ist klar, dass
\begin{equation}\label{l2}
\xi_{n-1}
< \xi_n+\xi_{n-1} = \xi (\pmb{z}_n-\pmb{z}_{n-1}).
\end{equation}
Aus  (\ref{l1},\ref{l2}) folgt, dass $\pmb{z}_{n-1}$ ein inner Punkt f\"ur $\Omega_1$ ist.

Wegen $ \frac{\xi_{n-1}}{\xi_n} = \alpha_{n+1}>
a_{n+1} \ge 2$ haben wir 
\begin{equation}\label{l10}
\xi_n < \xi_{n-1}-\xi_n = \xi (\pmb{z}_{n-1}+\pmb{z}_n)
,
\end{equation}
und es ist klar, dass
\begin{equation}\label{l20}
q_{n}< q_n +q_{n-1}.
\end{equation}
Aus  (\ref{l10},\ref{l20}) folgt, dass $\pmb{z}_{n}$ ein inner Punkt f\"ur $\Omega_2$ ist.

Nun es gen\"ugt zu beweisen, dass $\Pi$ keine Gitterpunkten hat, au\ss er
 $
 \pmb{z}_{n-1}, \pmb{z}_n
$
und
$\pmb{0}, 
\pmb{z}_{n} - \pmb{z}_{n-1}, \pmb{z}_n + \pmb{z}_{n-1}
$
auf dem Rand.

 Das Gitter $\mathbb{Z}^2$ teilt sich in den Linien
$$
\mathbb{Z}^2 =
\bigcup_{u \in \mathbb{Z}} \,\ell_u,\,\,\,\
\ell_u =\{
\pmb{z}\in \mathbb{Z}^2:
\,\,
\pmb{z} = u\pmb{z}_{n-1}+v(\pmb{z}_n -\pmb{z}_{n-1}), \,\, v \in \mathbb{Z}\}.
$$

 Falls $\pmb{z} \in \ell_u $ mit $|u|\ge 2$, die Ungleichung
 $$
 \xi (\pmb{z}) > 2\xi_{n-1}-\xi_n> \xi_{n-1}+\xi_n = \xi (\pmb{z}_{n-1}-
 \pmb{z}_n)$$
 gilt (die letzte Ungleichung aus
 $ \frac{\xi_{n-1}}{\xi_n} =\alpha_{n+1} > a_{n+1} \ge 2$
 folgt). Sodass f\"ur diese $u$  haben wir
 $ \ell_u \cap \Pi =\varnothing$.

 Es ist klar, dass
 $\ell_{-1}\cap\Pi =\{  
 \pmb{z}_n -\pmb{z}_{n-1}\}$,\,\,
  $\ell_{0}\cap\Pi =\{ \pmb{0}, 
 \pmb{z}_n \}$,\,\,
  $\ell_{1}\cap\Pi =\{  
 \pmb{z}_{n-1},\pmb{z}_n +\pmb{z}_{n-1}\}$.
 
 Daraus folgt
 $$
 \Pi \cap\mathbb{Z}^2 =
 \{  \pmb{0}, 
 \pmb{z}_{n-1},\pmb{z}_n, \pmb{z}_n -\pmb{z}_{n-1}, \pmb{z}_n +\pmb{z}_{n-1}\}
 .
 $$
 Damit ist Hilfssatz 2 bewiesen.$\Box$

  {\bf Hilfssatz 3.}
\,\,{\it
Sei $a_n \ge 2$ und $a_{n+1}=1$. Dann 
$$
q_{n}-q_{n-1},\,\,\,
2q_{n-1}+q_{n}
$$
sind zwei sukzessive Elemente der Sequenz $\hbox{\got Q}$.
}

Beweis  (siehe Fig. 3).\, Beweis des Hilfssatzes 3 verl\"auft analog zu dem Beweis des Hilfssatzes 2.
Betrachten wir die Parallelorgamme
$$
 \Pi = 
\{
(x,y) \in \mathbb{R}^2:\,\,
0\le x\le  2q_{n-1}+q_n,\,\,\,
|\alpha x-y|\le \xi_{n-1}+\xi_{n}\},
$$
$$
\Omega_1 = 
\{
(x,y) \in \mathbb{R}^2:\,\,
0\le x\le  q_{n}-q_{n-1},\,\,\,
|\alpha x-y|\le \xi_{n-1}+\xi_{n}\}\subset \Pi,
$$
$$
\Omega_2 = 
\{
(x,y) \in \mathbb{R}^2:\,\,
0\le x\le  2q_{n-1}+q_n,\,\,\,
|\alpha x-y|\le \xi_{n-1}-\xi_{n}\}\subset\Pi.
$$
Nun $\pmb{z}_{n-1}$ ist ein inner Punkt f\"ur $\Omega_1$ t, und
$\pmb{z}_n,\pmb{z}_{n+1}$  sind innere Punkte f\"ur $\Omega_2$.
F\"ur das Gitter haben wir
$$
\mathbb{Z}^2 =
\bigcup_{u \in \mathbb{Z}} \,\ell_u,\,\,\,\
\ell_u =\{
\pmb{z}\in \mathbb{Z}^2:
\,\,
\pmb{z} = u\pmb{z}_{n-1}+v\pmb{z}_{n}, \,\, v \in \mathbb{Z}\},
$$
und
$$
\ell_{-1}\cap \Pi = \{\pmb{z}_n -\pmb{z}_{n-1}\},
\ell_0\cap \Pi =\{ \pmb{z}_n\},
\ell_{1} \cap \Pi =
\{ \pmb{z}_{n-1}, \pmb{z}_{n+1}\},
\ell_{2}\cap \Pi = \{2\pmb{z}_{n-1} +\pmb{z}_{n}\},\,\,\,
\ell_u \cap \Pi =\varnothing, \,u \neq -1,0,1,2.
$$
Sodass
$$
\Pi \cap \mathbb{Z}^2 =
 \{\pmb{0}, \pmb{z}_{n-1},\pmb{z}_n,\pmb{z}_{n+1}, \pmb{z}_n-\pmb{z}_{n-1}, 2\pmb{z}_{n-1}+\pmb{z}_n\}.
$$
Hilfssatz 3 ist bewiesen.$\Box$

  {\bf Hilfssatz 4.}
\,\,{\it
Sei $a_n \ge 2$ und $a_{n+1}=a_{n+2}=1$. Dann 
$$
q_{n}-q_{n-1},\,\,\,
2q_{n-1}+q_{n},\,\,\,
2q_n,\,\,\, 2q_{n+1}
$$
sind vier sukzessive Elemente der Sequenz $\hbox{\got Q}$.
}

\begin{figure}[h]
  \centering
  \begin{tikzpicture}[scale=1.4]
    
    \draw[color=black] (0.4,2.6) -- (3.6,-1.4);
    \draw[color=black] (0,0) -- (4,1.2);
    \draw[color=black] (0,0) -- (1.6,-2);
 \draw[color=black] (4,1.2) -- (5.6,-0.8);
  \draw[color=black] (0.8,-1) -- (4.8,0.2);
  \draw[color=black] (1.6,-2) -- (5.6,-0.8);

 \node[fill=black,circle,inner sep=1.2pt]   at (0,0) {};
 \node[fill=black,circle,inner sep=1.2pt] at (0.4,2.6) {};
  \node[fill=black,circle,inner sep=1.2pt]   at (1.6,-2) {};
  \node[fill=black,circle,inner sep=1.2pt]   at (2,0.6) {};
  \node[fill=black,circle,inner sep=1.2pt] at (2.8,-0.4) {};
  \node[fill=black,circle,inner sep=1.2pt]   at (0.8,-1) {};
  \node[fill=black,circle,inner sep=1.2pt]   at (4.8,0.2) {};

  \node[draw=black,circle,inner sep=2pt] at (2,0.6) {};
  \node[draw=black,circle,inner sep=2pt] at (0.8,-1) {};
  \node[draw=black,circle,inner sep=2pt] at (2.8,-0.4) {};
    \node[draw=black,circle,inner sep=2pt] at (4.8,0.2) {};

  \draw[color=black] (0,0) -- (6.7,0);
         \draw (6.3,0.1)node[]{\begin{small}$y=\alpha x$\end{small}};
  
  \draw[color=black] (0,-2.1) -- (0,2.7);

    \draw[color=black] (0,1.4) -- (4,1.4);
    \draw[color=black] (0,-1.4) -- (4,-1.4);
    \draw[color=black] (4,1.4) -- (4,-1.4);

    \draw[color=black] (0,1.6) -- (3.6,1.6);
  \draw[color=black] (0,-1.6) -- (3.6,-1.6);
 \draw[color=black] (3.6,1.6) -- (3.6,-1.6);

  \draw[color=black] (0,1.2) -- (5.6,1.2);
  \draw[color=black] (0,-1.2) -- (5.6,-1.2);
 \draw[color=black] (5.6,1.2) -- (5.6,-1.2);

    \draw[color=black] (0,1.6) -- (1.2,1.6);
  \draw[color=black] (0,-1.6) -- (1.2,-1.6);

  \draw (1.9,1.75)node[]{\begin{small}$\pmb{z}_{n}-\pmb{z}_{n-1}$\end{small}};
  \draw (4.35,-1.6)node[]{\begin{small}$2\pmb{z}_{n-1}+\pmb{z}_{n}$\end{small}};
  \draw (4.3,1.4)node[]{\begin{small}$2\pmb{z}_{n}$\end{small}};
  \draw (6.1,-0.8)node[]{\begin{small}$2\pmb{z}_{n+1}$\end{small}};

 \draw (0.95,-0.75)node[]{\begin{small}$\pmb{z}_{n-1}$\end{small}};
 \draw (2.12,0.8)node[]{\begin{small}$\pmb{z}_{n}$\end{small}};
 \draw (3.27,-0.52)node[]{\begin{small}$\pmb{z}_{n+1}$\end{small}};
 \draw (5.2,0.2)node[]{\begin{small}$\pmb{z}_{n+2}$\end{small}};
 
  \node[draw=black,fill=red,circle,inner sep=2pt] at (1.2,1.6) {};
  \node[draw=black,fill=red,circle,inner sep=2pt] at (4,1.2) {};
   \node[draw=black,fill=red,circle,inner sep=2pt] at (5.6,-0.8) {};
    \node[draw=black,fill=red,circle,inner sep=2pt] at (3.6,-1.4) {};
  \end{tikzpicture}
  
   \caption{$a_n \ge 2,  a_{n+1}=a_{n+2}=1$} \label{fig4}
\end{figure}
Beweis  (siehe Fig. 3).\,
Aus Hilfssatz 3 folgt, dass 
$q_{n}-q_{n-1},2q_{n-1}+q_n \in \hbox{\got Q}$.
Sei
$$
\Omega
=
\{
\pmb{z} = (x,y) \in \mathbb{R}^2
:\,\,\,
0\le x\le 2q_n,\,\,
|x\alpha - y|\le
2\xi_{n-1}-\xi_n =\xi(2\pmb{z}_{n-1}+\pmb{z}_n)
\}.
$$
Es gen\"ugt  zu beweisen, dass
\begin{equation}\label{e2}
  \xi (2\pmb{z}_n) = 2\xi_n < 2\xi_{n-1} - \xi_n =
 \xi(2\pmb{z}_{n-1}+\pmb{z}_n).
\end{equation}
und
\begin{equation}\label{e1}
 \Omega \cap \mathbb{Z}^2
 =
 \{\pmb{0},
 \pmb{z}_{n-1},\pmb{z}_n,\pmb{z}_{n+1},2\pmb{z}_{n-1}+\pmb{z}_n,
 2\pmb{z}_n\},
\end{equation}
wo Punkte $\pmb{z}_{n+1},2\pmb{z}_{n-1}+\pmb{z}_n$ auf dem Rand liegen und
$\pmb{z}_{n-1},\pmb{z}_n,\pmb{z}_{n+1}\in {\rm int}\, \Omega$.
 
F\"ur die Gitterpunkte in $\Omega$ ist die Behauptung (\ref{e1}) klar.
Die Ungleichung (\ref{e2}) folgt aus
$
\frac{\xi_{n-1}}{\xi_n} =\alpha_{n+1} = [1;1,...] =
1+\frac{1}{1+\cdots} > \frac{3}{2}$.$\Box$

  {\bf Hilfssatz 5.}
\,\,{\it
Sei $a_n =1$ und $a_{n+1}\ge 2$. Dann 
$$
2q_{n-2}+q_{n-1},\,\,\,
q_{n-1}+q_{n}$$
sind zwei sukzessive Elemente der Sequenz $\hbox{\got Q}$.
}

\begin{figure}[h]
  \centering
  \begin{tikzpicture}[scale=2.2]

     \draw[color=black] (0,1) -- (1.7,1);
     \draw[color=black] (1.88,1) -- (2.2,1);
    \draw[color=black] (2.2,1) -- (2.2,-1);
      \draw (1.8,1.01)node[]{$\Pi$};
     
     
     \draw[color=black] (0,-1) -- (2.2,-1);

        \draw[color=black] (0,0) -- (3.45,0);
         \draw (3.1,0.1)node[]{\begin{small}$y=\alpha x$\end{small}};
        
        \draw[color=black] (0,-1.2) -- (0,1.2);

   \draw[color=black] (0,0) -- (1.2,0.1)  ;
  \draw[color=black] (-0.2,-0.9) -- (0.2,0.9)  ;
      \draw[color=black] (0.2,0.9) -- (1.4,1)  ;
       \draw[color=black] (-0.2,-0.9) -- (3.4,-0.6)  ;
  \draw[color=black] (1,-0.8) -- (1.4,1)  ;

  \node[fill=black,circle,inner sep=1.2pt] at (0,0) {};
  
  \node[fill=black,circle,inner sep=1.2pt] at (0.2,0.9) {};
  \node[fill=black,circle,inner sep=1.2pt] at (-0.2,-0.9) {};
  \node[fill=black,circle,inner sep=1.2pt] at (1.2,0.1) {};
  \node[fill=black,circle,inner sep=1.2pt] at (1,-0.8) {};
  \node[fill=black,circle,inner sep=1.2pt] at (3.4,-0.6) {};

   \node[draw=black,circle,inner sep=2pt] at (0.2,0.9) {};
  \node[draw=black,circle,inner sep=2pt] at (1,-0.8) {}; 
   \node[draw=black,circle,inner sep=2pt] at (1.2,0.1) {};

   \draw[dashed,color=black] (1.4,0.98) -- (1.4, -0.4);
   \draw[dashed,color=black] (1.4,-0.98) -- (1.4, -0.6); 
   \draw[dashed,color=black] (0,0.98) -- (1.4, 0.98);
   \draw[dashed,color=black] (0,-0.98) -- (1.4, - 0.98);
   \draw (1.4,-0.52)node[]{$\Omega_1$};

  \draw[dashed,color=black] (2.18,0.68) -- (2.18, - 0.68);
    
    \draw[dashed,color=black] (0,0.68) --  (1.67,0.68);
    \draw[dashed,color=black] (1.89,0.68) --  (2.18,0.68);
    \draw[dashed,color=black] (1,-0.68) --  (2.18,-0.68); 
    \draw[dashed,color=black] (0,-0.68) --  (0.6,-0.68);   
     \draw (1.8,0.66)node[]{$\Omega_2$};

     \draw (2.65,-0.8)node[]{\begin{small}$\pmb{z}_{n-1}+\pmb{z}_n$\end{small}};
     
     \draw (1.2,1.15)node[]{\begin{small}$2\pmb{z}_{n-2}+\pmb{z}_{n-1}$\end{small}};
   
     \draw (1.08,0.2)node[]{\begin{small}$\pmb{z}_n$\end{small}};
    
   \draw (0.8,-0.7)node[]{\begin{small}$\pmb{z}_{n-1}$\end{small}};
 
   \draw (0.4,0.79)node[]{\begin{small}$\pmb{z}_{n-2}$\end{small}};
 
    \node[draw=black,fill=red,circle,inner sep=2.2pt] at (1.4,1) {};

   \node[draw=black,fill=red,circle,inner sep=2.2pt] at (2.2,-0.7) {};

  \end{tikzpicture}
  
   \caption{$a_n = 1,a_{n+1}\ge 2$} \label{fig100}
\end{figure}

Der Beweis verl\"auft analog  (siehe Fig. 4).
Wir bertachten 
die Punkte $2\pmb{z}_{n-2}+\pmb{z}_{n-1}$ und $ \pmb{z}_{n-1}+\pmb{z}_n$ auf dem Rand des Parellogrammes
$$
\Pi =
\{
(x,y) \in \mathbb{R}^2:\,\,
0\le x\le  q_{n-1}+q_n,\,\,\,
|\alpha x-y|\le 2\xi_{n-2}-\xi_{n-1}\},
$$
und beweisen, dass
$
\Pi \cap \mathbb{Z}^2 =
\{\pmb{0}, \pmb{z}_{n-2},\pmb{z}_{n-1},\pmb{z}_{n}, 2\pmb{z}_{n-2}+\pmb{z}_{n-1}, \pmb{z}_{n-1}+\pmb{z}_n\}$.$\Box$

   {\bf Hilfssatz 6.}
\,\,{\it
Sei $a_{n-1}=a_n =1$ und $a_{n+1}\ge 2$. Dann 
$$
2q_{n-2},\,\,\,
2q_{n-1},\,\,\,
2q_{n-2}+q_{n-1},
\,\,\,
q_{n-1}+q_{n}
$$
sind vier sukzessive Elemente der Sequenz $\hbox{\got Q}$.
}

 Beweis analog zum Beweis des Hilfssatzes 4 ist.$\Box$

    {\bf Hilfssatz 7.}
\,\,{\it
Sei $a_n =a_{n+1}=a_{n+2} =1$. Dann 
$$2q_{n},\,\,\,
2q_{n+1}$$
sind zwei sukzessive Elemente der Sequenz $\hbox{\got Q}$.
}
\begin{figure}[h]
  \centering
  \begin{tikzpicture}[scale=1.5]

     \draw[color=black] (0,1.2) -- (3.8,1.2);
     \draw[color=black] (4.15,1.2) -- (4.4,1.2);
      \draw (4,1.24)node[]{$\Pi$};
     
     \draw[color=black] (4.4,1.2) -- (4.4,-1.2);
     
     \draw[color=black] (0,-1.2) -- (0.15,-1.2);
     \draw[color=black] (1,-1.2) -- (4.4,-1.2);

        \draw[color=black] (0,0) -- (5.7,0);
         \draw (5.3,0.1)node[]{\begin{small}$y=\alpha x$\end{small}};
        
        \draw[color=black] (0,-2.1) -- (0,1.8);

       \node[draw=black,circle,inner sep=2pt] at (0.6,1.6) {};

    \draw[color=black] (0.6,1.6) -- (1.4,0.6)  ;

    \draw (0.8,1.8)node[]{\begin{small}$\pmb{z}_{n-2}$\end{small}};

    \draw[color=black] (0,0) -- (1.4,0.6) -- (2.8,1.2)  ;
  \node[fill=black,circle,inner sep=1.2pt] at (0,0) {};

    \draw (2.9,1.4)node[]{\begin{small}$2\pmb{z}_{n}$\end{small}};
   
    \draw (1.45,0.8)node[]{\begin{small}$\pmb{z}_{n}$\end{small}}; 
   
    \draw[color=black] (0,0) -- (0.8,-1) -- (1.6,-2)  ;
  \node[fill=black,circle,inner sep=1.2pt] at (0.8,-1) {};
    \draw (0.5,-1.1)node[]{\begin{small}$\pmb{z}_{n-1}$\end{small}};
   \node[draw=black,circle,inner sep=2pt] at (0.8,-1) {};

   \node[fill=black,circle,inner sep=1.2pt] at (1.6,-2) {};
    \draw[color=black] (0,-1.2) -- (4.4,-1.2);
      \draw (2.35,-0.14)node[]{\begin{small}$\pmb{z}_{n+1}$\end{small}};
   
    \draw[color=black] (0.6,1.6) -- (1.4,0.6) -- (2.2,-0.4) -- (3, -1.4) ;
    \node[fill=black,circle,inner sep=1.2pt] at (0.6,1.6) {};
     \node[fill=black,circle,inner sep=1.2pt] at (1.4,0.6) {};
       \node[fill=black,circle,inner sep=1.2pt] at (2.2,-0.4) {};
          \node[fill=black,circle,inner sep=1.2pt] at (3,-1.4) {};
      \node[draw=black,circle,inner sep=2pt] at (1.4,0.6) {};
       \node[draw=black,circle,inner sep=2pt] at (2.2,-0.4) {};

    \draw[color=black] (0.8,-1) -- (2.2,-0.4) -- (3.6,0.2)  ;

    \draw[color=black] (2.8,1.2) -- (3.6,0.2) -- (4.4,-0.8)  ;
  
    \draw[color=black] (1.6,-2) -- (3,-1.4) -- (4.4,-0.8)  ;

     \node[fill=black,circle,inner sep=1.2pt] at (3.6,0.2) {};
    \draw (4,0.2)node[]{\begin{small}$\pmb{z}_{n+2}$\end{small}};
   \node[draw=black,circle,inner sep=2pt] at (3.6,0.2) {};

   \draw[dashed,color=black] (2.8,1.2) -- (2.8, 0.65);
   \draw[dashed,color=black] (2.8, 0.3) -- (2.8, - 1.2);
   \draw[dashed,color=black] (0,1.17) -- (2.8, 1.17);
   \draw[dashed,color=black] (0,-1.17) -- (2.8, - 1.17);
   \draw (2.8,0.5)node[]{$\Omega_1$};

  \draw[dashed,color=black] (4.37,0.8) -- (4.37, - 0.8);
    
    \draw[dashed,color=black] (0,0.8) --  (1.3,0.8);
   \draw[dashed,color=black] (1.6,0.8) --  (3.6,0.8);
     \draw[dashed,color=black]  (4,0.8) --   (4.37, 0.8);
   \draw[dashed,color=black] (0,-0.8) -- (4.37, - 0.8);
      \draw (3.8,0.8)node[]{$\Omega_2$};
  
   \draw (4.85,-0.8)node[]{\begin{small}$2\pmb{z}_{n+1}$\end{small}};
  
    \node[draw=black,fill=red,circle,inner sep=2.2pt] at (2.8,1.2) {};

   \node[draw=black,fill=red,circle,inner sep=2.2pt] at (4.4,-0.8) {};
 
  \end{tikzpicture}
  
   \caption{$a_n = a_{n+1}=a_{n+2}=1$} \label{fig7}
\end{figure}
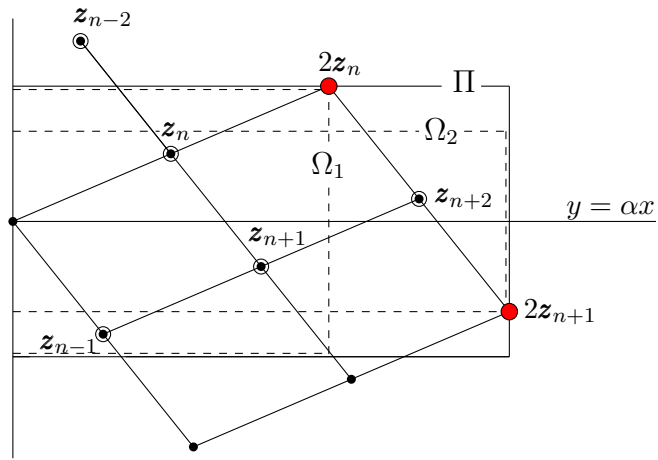

Beweis  (siehe Fig. 5).\,
Betrachten wir die Parallelogramme
$\Pi =\Pi [2\pmb{z}_n, 
2\pmb{z}_{n+1}], \Omega_1 =\Omega[2\pmb{z}_n],\Omega_2= \Omega[2\pmb{z}_{n+1}]$
(Fig. \ref{fig7}).

Es ist klar, dass $2\xi_{n+1} < 2\xi_n$,
und die Punkte
$2\pmb{z}_n,
2\pmb{z}_{n+1}$
auf dem Rand des Parallelogrammes $\Pi$ liegen.

Aus $\frac{\xi_n}{\xi_{n+1}}=\alpha_{n+2}<a_{n+2}+1 = 2$
folgt $\xi_n <2\xi_{n+1}$ und
$$
\xi_{n+2}< \xi_{n+1}< \xi_{n} <2\xi_{n+1}.$$
Au\ss erdem
$$
q_n <q_{n+1} < q_{n+2} <2q_{n+1}.
$$
Sodass $\pmb{z}_n,\pmb{z}_{n+1}, \pmb{z}_{n+2}$
sind innere Punkte des Parallelogrammes $\Omega_2$.

Aus $\frac{\xi_{n-1}}{\xi_n}=\alpha_{n+1}<a_{n+1}+1 = 2$ und
$q_{n+1} < 2q_n <q_{n+2}$ folgt, dass
$\pmb{z}_{n-1},\pmb{z}_{n}, \pmb{z}_{n+1}$
sind innere Punkte des Parallelogrammes $\Omega_1$.

Wegen $\frac{\xi_{n-1}}{\xi_n} =\alpha_{n+1} = 1 +\frac{1}{1+\cdots} > \frac{3}{2}$
haben wir
$$
\xi (2\pmb{z}_{n-1}+\pmb{z}_n) =
2\xi_{n-1} - \xi_n > 2\xi_n =
\xi (2\pmb{z}_n).
$$
Au\ss erdem 
$$
\xi (\pmb{z}_{n-2})
=
\xi_{n-2}> 2\
\xi_n = \xi (2\pmb{z}_n).
$$
Wir betrachten die Teilung 
$
\mathbb{Z}^2 =
\bigcup_{u \in \mathbb{Z}} \,\ell_u,\,\,\,\
\ell_u =\{
\pmb{z}\in \mathbb{Z}^2:
\,\,
\pmb{z} = u\pmb{z}_{n-1}+v\pmb{z}_{n}, \,\, v \in \mathbb{Z}\},
$ und ziehen
$$
\Pi \cap \mathbb{Z}^2 =
 \{\pmb{0}, \pmb{z}_{n-1},\pmb{z}_n,\pmb{z}_{n+1}, 2\pmb{z}_n, \pmb{z}_{n+2},2\pmb{z}_{n+1}\}.
$$
Hilfssatz 7  bewiesen ist.$\Box$

{\bf 3.   \"Uber die Funktion  $\psi_\alpha^{[2]*} (t)$.}

 Es gibt eine einfache
Regel, die die Sequenz $ \hbox{\got X}$ aus der Sequenz $\hbox{\got A}$ konstruirt.

 {\bf Satz 2.}
 \,\,
 {\it Die Sequenz $\hbox{\got X}$ wird mit dem folgenden Regel  aus der Sequenz $\hbox{\got A}
 $ erhalten:

  \noindent
 1. jedes Element $ a_n \ge 2$   wird durch 
 $a_n-1$ sukzessive Zahlen $q_{n-2}+jq_{n-1},\,\,\,
1\le j \le a_{n}-1$ ersetzen;
 
\noindent
 2. jedes Element $ a_n =1$   wird  durch 
 eine Zahl $2q_{n-2}+q_{n-1}$ ersetzen.}

 {\bf Bemerkung 2.} \, {\it Die Regel gibt uns  verschiedene Elemente der 
 Sequenz $\hbox{\got X}$. Die einzige Ausnahme ist im Fall
 $a_1=a_2=1, a_3 \ge 2$, als im Punkt 2 der Bemerkung 1.}

Wir nehmen die Notation des Haupthilfssatz 1 an.
Aus der Definition der Funktion $\psi_\alpha^{[2]*}$ folgt:

{\bf Haupthilfssatz 2.}\,\,{\it
Um zu zeigen, dass 
$x_1< x_2$  zwei sukzessive Elemente aus $\hbox{\got Q}$ sind,
genügt es, folgendes zu beweisen:

\noindent
$\bullet$ f\"ur einige  $y_1,y_2 \in \mathbb{Z}$ die Gitterpunkte
$\pmb{w}_1 = (x_1,y_1), 
\pmb{w}_2 = (x_2,y_2)
$
auf dem Rand des Parallelogrammes
$ 
\Pi [\pmb{w}_1, 
\pmb{w}_2].
$
liegen;

\noindent
$\bullet$
das Parallelorgamm 
$
\Omega [\pmb{w}_1]  $
hat einen inneren Punkt $\pmb{z}_{\nu_1}  = (p_{\nu_1}, q_{\nu_1})$;
 
 \noindent
$\bullet$
das Parallelorgamm 
$
\Omega [
\pmb{w}_2] 
$
hat einen inneren Punkt $\pmb{z}_{\nu_2} = (p_{\nu_2}, q_{\nu_2})$ auch;
 
\noindent
$\bullet$
falls $\pmb{z} \in \Pi [\pmb{w}_1, 
\pmb{w}_2] $
ein Gitterpunt ist, dann
$\pmb{z} \in \{\pmb{0}, \pmb{w}_1,\pmb{w}_2\}$ oder
$\pmb{z} = k\cdot \pmb{z}_\nu$ mit einem $\nu\ge -1$ und $k \in \mathbb{Z}_+$.}

Die folgende Hilfss\"atze sind klar. Hier geben wir keine Beweise.
Elemente der Sequentz $\hbox{\got X}$ entsprechen gr\"unen Punkten auf den Figuren 6 und 7.

    {\bf Hilfssatz 8.}\footnote{Der Autor dankt Prof. V. Bykovskii, der ihm  diese Behauptung   erkl\"art hatte.}
\,\,{\it
Sei $a_n \ge 2$. Dann 
$$q_{n-2}+jq_{n-1},\,\,\,
1\le j \le a_{n}-1$$
sind $a_{n}-2$ sukzessive Elemente der Sequenz $\hbox{\got X}$ (siehe Fig. 6).
}

\begin{figure}[h]
  \centering
  \begin{tikzpicture}[scale=1.6]

  \draw[color=black] (0,0) -- (7.7,0);
         \draw (7.35,-0.2)node[]{\begin{small}$y=\alpha x$\end{small}};
  
  \draw[color=black] (0,-2.1) -- (0,1.7);

    \draw[color=black] (-0.4,-1.4) -- (0.4,1.4);
    \draw[color=black] (0.4,1.4) -- (8.2,-0.03);
 \draw[color=black] (0,0) -- (4.2,-0.77);
  \draw[color=black] (-0.4,-1.4) -- (2.6,-1.94);
 
    \draw[color=black] (0.8,-1.62) -- (1.6,1.18);
  
    \draw[color=black] (2,-1.84) -- (2.8,0.96);

   \node[fill=black,circle,inner sep=1.2pt]   at (7.6,0.08) {};
   
    \node[fill=black,circle,inner sep=1.2pt]   at (0,0) {};
     \node[fill=black,circle,inner sep=1.2pt]   at (-0.4,-1.4) {};
    \node[fill=black,circle,inner sep=1.2pt]   at (0.4,1.4) {};   
 
         \node[fill=black,circle,inner sep=1.2pt]   at (0.8,-1.62) {};
           \node[fill=black,circle,inner sep=1.2pt]   at (2,-1.84) {};

      \node[fill=black,circle,inner sep=1.2pt]   at (1.2,-0.22) {};

      \node[fill=black,circle,inner sep=1.2pt]   at (3.6,-0.66) {};

  \node[draw=black,circle,inner sep=2pt] at (0.4,1.4) {};
  \draw (0.75,1.5)node[]{$\pmb{z}_{n-2}$};

  \node[draw=black,circle,inner sep=2pt] at (7.6,0.08) {};
   \draw (7.8,0.2)node[]{$\pmb{z}_{n}$};

  \node[draw=black,circle,inner sep=2pt] at (1.2,-0.22) {};
   \draw (1.55,-0.15)node[]{$\pmb{z}_{n-1}$};

    \draw[color=black] (0,1.18) -- (2.8,1.18);
    \draw[color=black] (0,-1.18) -- (1.7,-1.18);
    \draw[color=black] (2.1,-1.18) -- (2.8,-1.18);
    \draw[color=black] (2.8,1.18) -- (2.8,-0.65);
    \draw[color=black] (2.8,-0.74) -- (2.8,-1.18);
    
        \draw (1.9,-1.18)node[]{$\Pi_1$};

    \draw[color=black] (0,0.96) -- (4,0.96);
    \draw[color=black] (0,-0.96) -- (4,-0.96);
     
    \draw[color=black] (4,0.96) -- (4,-0.96);

    \draw[color=black] (0,0.74) -- (5.2,0.74);
    \draw[color=black] (0,-0.74) -- (5.2,-0.74);
     
    \draw[color=black] (5.2,0.74) -- (5.2,-0.74);

       \draw[color=black] (0,0.52) -- (6.4,0.52);
    \draw[color=black] (0,-0.52) -- (6.4,-0.52);
    
    \draw[color=black] (6.4,0.52) -- (6.4,-0.52);

     \draw (2.35,1.35)node[]{$\pmb{z}_{n-2}+\pmb{z}_{n-1}$};  

          \draw (2.8,-0.65)node[]{$2\pmb{z}_{n-1}$};  
    
         \draw (6.97,0.44)node[]{$\pmb{z}_{n}-\pmb{z}_{n-1}$};  
    
  \node[draw=black,fill=green,circle,inner sep=2pt] at (1.6,1.18) {};
   \node[draw=black,fill=green,circle,inner sep=2pt]   at (2.8,0.96) {};

     \node[draw=black,fill=green,circle,inner sep=2pt] at (4,0.74) {};
   \node[draw=black,fill=green,circle,inner sep=2pt]    at (5.2,0.52) {};
   
   \node[draw=black,fill=green,circle,inner sep=2pt] at (6.4,0.3) {};

  \node[draw=black,fill=red,circle,inner sep=2pt] at (2.4,-0.44) {};

  \end{tikzpicture}
  
   \caption{$a_n =6$} \label{fig6}
\end{figure}

    {\bf Hilfssatz 9.}
\,\,{\it
Sei $a_n, a_{n+1} \ge 2$. Dann 
$$q_{n}-q_{n-1},\,\,\,
q_n+q_{n-1}$$
sind zwei sukzessive Elemente der Sequenz $\hbox{\got X}$ (siehe Fig 2.).
}

    {\bf Hilfssatz 10.}
\,\,{\it
Sei $a_n \ge 2, a_{n+1}=1$. Dann 
$$q_{n}-q_{n-1},\,\,\,
2q_{n-1}+q_n$$
sind zwei sukzessive Elemente der Sequenz $\hbox{\got X}$ (siehe Fig 3.).
}

   {\bf Hilfssatz 11.}
\,\,{\it
Sei $a_n =1, a_{n+1}\ge 2$. Dann 
$$2q_{n-2}-q_{n-1},\,\,\,
q_{n+1}-q_n$$
sind zwei sukzessive Elemente der Sequenz $\hbox{\got X}$ (siehe Fig 4.).
}

    {\bf Hilfssatz 12.}
\,\,{\it
Sei $a_n =a_{n+1}=1$. Dann 
$$2q_{n-2}+q_{n-1},\,\,\,
2q_{n-1}+q_n$$
sind zwei sukzessive Elemente der Sequenz $\hbox{\got X}$
(seihe Fig 7.).}

\begin{figure}[h]
  \centering
  \begin{tikzpicture}[scale=1.6]

     \draw[color=black] (0,1.37) -- (3.8,1.37);
     \draw[color=black] (4.15,1.37) -- (5,1.37);
      \draw (4,1.36)node[]{$\Pi$};
     
     \draw[color=black] (5,1.37) -- (5,-1.37);
     
     \draw[color=black] (0,-1.37) -- (5,-1.37);

        \draw[color=black] (0,0) -- (6.45,0);
         \draw (6.1,0.1)node[]{\begin{small}$y=\alpha x$\end{small}};
        
        \draw[color=black] (0,-2.1) -- (0,1.8);

    \node[fill=black,circle,inner sep=1.2pt] at (0.6,1.63) {};
                       

   \draw[color=black] (0,0) -- (1.4,0.63) -- (2.8,1.26) -- (4.2,1.89)  ;

    \draw[color=black] (5.8,-0.11) -- (4.2,1.89)  ;

  \node[fill=black,circle,inner sep=1.2pt] at (0,0) {};
  
  \node[fill=black,circle,inner sep=1.2pt] at (4.2,1.89) {};

    \draw (2.85,1.5)node[]{\begin{small}$2\pmb{z}_{n-1}$\end{small}};
   
    \draw (1.78,0.53)node[]{\begin{small}$\pmb{z}_{n-1}$\end{small}}; 
   
    \draw[color=black] (0,0) -- (0.8,-1) -- (1.6,-2)  ;
  \node[fill=black,circle,inner sep=1.2pt] at (0.8,-1) {};
    \draw (0.6,-1.15)node[]{\begin{small}$\pmb{z}_{n-2}$\end{small}};
   \node[draw=black,circle,inner sep=2pt] at (0.8,-1) {};

   \node[fill=black,circle,inner sep=1.2pt] at (1.6,-2) {};
    
      \draw (2.25,-0.2)node[]{\begin{small}$\pmb{z}_{n}$\end{small}};
   
    \draw[color=black] (0.6,1.63) -- (1.4,0.63) -- (2.2,-0.37) -- (3, -1.37) ;
     \node[fill=black,circle,inner sep=1.2pt] at (1.4,0.63) {};
       \node[fill=black,circle,inner sep=1.2pt] at (2.2,-0.37) {};
          \node[fill=black,circle,inner sep=1.2pt] at (3,-1.37) {};
      \node[draw=black,circle,inner sep=2pt] at (1.4,0.63) {};
       \node[draw=black,circle,inner sep=2pt] at (2.2,-0.37) {};

    \draw[color=black] (0.8,-1) -- (2.2,-0.37) -- (3.6,0.26) -- (5,0.89) ;

    \draw[color=black] (2.8,1.26) -- (3.6,0.26) -- (4.4,-0.74)  ;
  
    \draw[color=black] (1.6,-2) -- (3,-1.37) -- (4.4,-0.74) -- (5.8,-0.11) ;

   \node[fill=black,circle,inner sep=1.2pt] at (5.8,-0.11) {};

     \node[fill=black,circle,inner sep=1.2pt] at (3.6,0.26) {};
    \draw (4,0.2)node[]{\begin{small}$\pmb{z}_{n+1}$\end{small}};
   \node[draw=black,circle,inner sep=2pt] at (3.6,0.26) {};

   \draw[dashed,color=black] (3,1.34) -- (3, 0.6);
   
   \draw[dashed,color=black] (3,0.4) -- (3, -1.34);
   
   \draw[dashed,color=black] (0,1.34) -- (3, 1.34);
   \draw[dashed,color=black] (0,-1.34) -- (3, - 1.34);
   \draw (3,0.5)node[]{$\Omega_1$};

  \draw[dashed,color=black] (4.97,0.89) -- (4.97, - 0.89);
    
    \draw[dashed,color=black] (0,0.89) --  (3.6,0.89);
   
     \draw[dashed,color=black]  (4,0.89) --   (4.97, 0.89);
   \draw[dashed,color=black] (0,-0.89) -- (4.97, - 0.89);
      \draw (3.8,0.89)node[]{$\Omega_2$};
  
   \draw (4.45,-0.48)node[]{\begin{small}$2\pmb{z}_{n}$\end{small}};
  
    \node[draw=black,fill=red,circle,inner sep=2.2pt] at (2.8,1.26) {};

   \node[draw=black,fill=red,circle,inner sep=2.2pt] at (4.4,-0.74) {};

    \node[draw=black,fill=green,circle,inner sep=2.2pt] at (5,0.89) {};

   \node[draw=black,fill=green,circle,inner sep=2.2pt] at (3,-1.37) {};
   \draw (3.78,-1.49)node[]{\begin{small}$2\pmb{z}_{n-2}+\pmb{z}_{n-1}$\end{small}}; 
   \draw (5.68,1)node[]{\begin{small}$2\pmb{z}_{n-1}+\pmb{z}_{n}$\end{small}}; 
  \end{tikzpicture}
  
   \caption{$a_n = a_{n+1}=1$} \label{fig71}
\end{figure}
 
Satz 2 folgt nach Hilfss\"atze 8 - 12.

 {\bf 4. \"Uber diophantische Spektra.}

 {\bf 4.1.  Asudr\"ucke mit Kettenbr\"uchen.}

Wir setzen
\begin{equation}\label{kk1}
\varkappa^1_n =
\varkappa^1_n (\alpha) =
(q_{n-2}+q_{n-1})(\xi_{n-2}-\xi_{n-1})
=
\frac{(1+\alpha_{n-1}^*)(\alpha_{n}-1)}
{\alpha_{n}+\alpha_{n-1}^*},
\end{equation}
\begin{equation}\label{kk2}
\varkappa^2_n =
\varkappa^2_n (\alpha) =
(q_n-q_{n-1})(\xi_{n-1}+\xi_n) =
\frac{(1-\alpha_n^*)(\alpha_{n+1}+1)}
{\alpha_{n+1}+\alpha_n^*}
,
\end{equation}
\begin{equation}\label{kk3}
\varkappa^3_n =
\varkappa^3_n (\alpha) =
(2q_{n-2}+q_{n-1})(2\xi_{n-2}-\xi_{n-1})=
\frac{(2\alpha_{n-1}^*+1)(2\alpha_{n}-1)}
{\alpha_{n}+\alpha_{n-1}^*}
,
\end{equation}
\begin{equation}\label{kk4}
\varkappa^4_n = \varkappa^4_n (\alpha) = 4q_{n-1}\xi_{n-1} = \frac{4}{\alpha_{n}+\alpha_{n-1}^*}.
\end{equation}

 Die Gleichungen  (\ref{kk1}) ---
(\ref{kk4}) folgen aus (\ref{2p}) und (\ref{2s}).
 Die letzte Gleichung in (\ref{kk4}) sehr bekanntlich ist (\cite{Cus}, Appendix 1), sie hat viele Anwendungen. 
  Nach dem Satz 1  folgt 
 $$
 \hbox{\got k}_\nu = \varkappa^j_n
 $$
 mit
einem $ n = n (\nu)$ und einem $ j \in \{1,2,3,4\}$.

 Beide Funktionen
 $
 \frac{(1+x)(y-1)}{x+y}$ und
 $
 \frac{(2x+1)(2y-1)}{x+y}
 $
  steigen 
  in 
   $x$ und $y$ im Bereich $x>0, y>1$.
  Die Funktion
 $
 \frac{(1-x)(y+1)}
{x+y}
 $
 f\"allt in $x$ und $y$.
  Daraus folgt
 \begin{equation}
  \label{man1}
  1\le \varkappa_n^3 \le 6,
 \end{equation}
 und falls  $a_n \ge 2$, haben wir
 \begin{equation}
  \label{man2}
  \frac{1}{2}\le \frac{a_{n}-1}{a_n}\le \varkappa_n^1, \varkappa_{n}^2 \le
  \frac{2a_n}{a_n+2}  
 < 2.
 \end{equation}

{\bf Bemerkung 3.} Falls $ a_n \ge 2$,
die Werte  $\varkappa_n^1 $ und $ \varkappa_n^2$ eine Symmetrie haben:
die Auswechselung
$
\alpha_n\mapsto \frac{1}{\alpha_{n-1}^*}, \alpha_{n-1}\mapsto \frac{1}{\alpha_n}
$
gibt uns
$
\varkappa_{n}^1\mapsto \varkappa_{n-1}^2.
$

{\bf Bemerkung 4.} 
Falls
$a_n =1$,
ist
\begin{equation}
 \label{oon}
 2q_{n-1}+q_{n} = q_{n-1}+q_{n+1}
\end{equation}
und
k\"onnen wir $\varkappa_n^3$ in einer symmetrischen Weise schreiben:
\begin{equation}\label{kk3k}
\varkappa^3_n =
(q_{n-2}+q_{n})(\xi_{n-2}-\xi_{n})
.
\end{equation}

 {\bf Bemerkung 5.} 
Falls
$a_n =2$,
haben wir 
\begin{equation}
 \label{twa}
 q_{n}-q_{n-1} = q_{n-1}+q_{n-2}
\end{equation}
und $ \varkappa^1_n =\varkappa^2_n$.

Sei 
$$
\hbox{\got l}_n^*
=
\hbox{\got l}_n^*(\alpha)
=
 \begin{cases}
  \varkappa_n^3\,\,\,\,\,\,\,\,\,\,\,\,\,\,\,\,\,\,\,\,\,\,\,\,\,\, 
  \text{falls}\,\,\, a_n =1,\cr
  \min (\varkappa_{n}^1,\varkappa_n^2) \,\,\, \text{falls} \,\,\,
 a_n \ge 2.
 \end{cases}
 $$
 Wegen
 $\min_{1\le j 
 \le a_{n}-1} (q_{n-2}+jq_{n-1})(\xi_{n-2}-j\xi_{n-1}) = 
  \min (\varkappa_{n}^1,\varkappa_n^2)$ 
 nach Satz 2 folgt
$$
 \hbox{\got j}^*(\alpha )=
\inf
\,\hbox{\got l}_n^*,\,\,\,\,
\hbox{\got k}^*=
\liminf_{n\to \infty}\,\hbox{\got l}_n^*.
$$

 {\bf Hilfssatz 13.}

 \noindent
   {\it 1. Es sei $ a_n =1, a_{n+1}\ge 2$ und
  $\varkappa_n^3 < \varkappa_{n+1}^1$. Dann ist $ a_{n+1}\ge 4  $ 
  und $\varkappa_n^3 >\varkappa_{n+1}^4$. 
  }

   \noindent
 {\it 2. Es sei  $ a_n =1, a_{n-1}\ge 2$ und
  $\varkappa_n^3 < \varkappa_{n-1}^2$. Dann ist $ a_{n-1}\ge 4  $ 
  und $\varkappa_n^3 >\varkappa_{n-1}^4$. 
  }
 
  \noindent   
{\it 3.    Es sei  $\alpha \not\sim\frac{1+\sqrt{5}}{2}$.
Dann ist 
   $$ \hbox{\got k}(\alpha ) =\liminf_{n\to \infty : a_n  \ge 2} 
   \,\min (
   \varkappa^1_n, \varkappa^2_n, \varkappa^4_n).$$
   }

    \noindent   
{\it
4. Falls $a_n \ge 8$, ist
 $$
 \varkappa_n^4 =
 \min_{1\le j\le 4}   \varkappa_n^j
 $$
 und f\"ur $\alpha$ mit unendlich vielen 
 $a_n \ge 8$ man hat
 \begin{equation}\label{lauf}
 \hbox{\got k}  (\alpha) =
 \liminf_{n\to \infty, \, a_n \ge 8}
 \varkappa_n^4 = 4\lambda (\alpha ).
 \end{equation}
   
   }

 Beweis. \,\,
 Hier werden wir nur die erste Behauptung beweisen. Beweis der zweiten Behauptung  verl\"auft analog.
  Behauptung 3 folgt aus  den Behauptungen 1,2.
Behauptung 4 ist klar. 
 
 Bemerken wir, dass $
\alpha_n = 1 + \frac{1}{\alpha_{n+1}}$ und $\frac{1}{\alpha_n^*} = 1 +\alpha_{n-1}^*$.
  Aus Definitionen (\ref{kk1},\ref{kk3})
 folgt
 $$
 (2q_{n-2}+q_{n-1}) (2\xi_{n-2}-\xi_{n-1}) < (q_{n-1}+q_n)(\xi_{n-1}-\xi_n),
 $$
sodass (\ref{2p},\ref{2s}) geben uns 
$$
(2\alpha_{n-1}^*+1)(2\alpha_n-1) < \left(1+\frac{1}{\alpha_n^*}\right) \left(1-\frac{1}{\alpha_{n+1}}
\right)
$$
und
$$
\frac{2+\alpha_{n+1}}{\alpha_{n+1}-1}
=
\frac{2\left(1+\frac{1}{\alpha_{n+1}}\right) - 1}{1-\frac{1}{\alpha_{n+1}}} =
\frac{2\alpha_n-1}{1-\frac{1}{\alpha_{n+1}}}
<
\frac{1+\frac{1}{\alpha_n^*}}{2\alpha_{n-1}^* +1} =
\frac{1+\frac{1}{1+\alpha_{n-1}^*}}{2\alpha_{n-1}^*+1}
<
\sup_{0<x<1}
\frac{2+x}{2x+1} = 2.
$$
 Die Ungleichung $ \alpha_{n+1} > 4$ folgt daraus und so $ a_{n+1} \ge 4$.
 Nun aus (\ref{man1}) und Definition (\ref{kk4}) folgt
 $
 \varkappa_{n+1}^4 
 < 1< \varkappa_n^3$.
 $\Box$

 {\bf 4.2. Beweis der Formeln (\ref{impo1}) und (\ref{impo2}).}

 Es ist klar dass f\"ur $\alpha = \sqrt{2}$ gilt
 $\lambda (\sqrt{2}) = \frac{1}{\sqrt{8}}$ und
 $
 \hbox{\got k} (\sqrt{2}) = \hbox{\got k}^* (\sqrt{2})  
 = \lim_{n\to \infty} \varkappa_n^2 = \frac{[1.\overline{2}]}{2} =\frac{1}{\sqrt{2}}$,
 und 
 (\ref{impo2})
 ist bewiesen.
 
 Um zu zeigen (\ref{impo1}), wir betrachten drei F\"alle.

 \noindent {\bf Fall 1$^0$.}\,
 Es gibt unendlich viele $n$ mit $ a_n \ge 4$. Dann $ \lambda(\alpha ) \le \frac{1}{4}$.
 Aus (\ref{man2}) folgt, dass $\varkappa^{1}_n, 
 \varkappa^{2}_n
\ge \frac{1}{2}$. Sodass nach Hilfssatz 13 (Behauptung 3)  haben wir $\hbox{\got k} (\alpha) \ge 2 \lambda (\alpha)$.

 \noindent {\bf Fall 2$^0$.}\,
 F\"ur alle hinreichend gro\ss en $n$, alle Teilenner $a_n$ sind  $\le 3$.
 Dann 
 $$\lambda (\alpha ) \le \frac{1}{3+2\cdot [0;\overline{3,1}]}
 <
 \frac{[\overline{1;3}]}{4} \le 
 \frac{1}{2} \cdot \min 
\left( \liminf_{n\to\infty : a_n \ge 2}\varkappa^{1}_n, 
 \liminf_{n\to\infty : a_n \ge 2}
 \varkappa^{2}_n \right) \le \frac{\hbox{\got k} (\alpha)}{2}
 $$
 und alles ist bewiesen.

  \noindent {\bf Fall 3$^0$.}\,
 F\"ur alle hinreichend gro\ss en $n$,
 man hat $ a_n \in \{1,2\}$. 
 F\"ur $\alpha \sim \frac{1+\sqrt{5}}{2}$  haben wir
 $\hbox{\got k} (\alpha) = 4 \lambda (\alpha)$. 
    F\"ur $\alpha \sim \sqrt{2}$  haben wir
 $\hbox{\got k} (\alpha) = 2 \lambda (\alpha)$. 
  Falls $\alpha \not\sim \frac{1+\sqrt{5}}{2}$ und $\alpha\not\sim \sqrt{2}$,
 nach der Struktur des Lagrengeschen Spektrums (Kapitel II aus \cite{Cas} oder \cite{Cus})
  haben wir
 $ \lambda (\alpha ) \le \frac{5}{\sqrt{221}} = 0.3363^+$
 und
 $ \liminf_{n\to\infty : a_n =2} \varkappa^1_n \ge \frac{[\overline{1;2}]}{2} = 0.683^+$.
 Sodass 
 $\hbox{\got k} (\alpha) \ge 2 \lambda (\alpha)$
 und alles ist bewiesen.

 {\bf 4.3. Das Spekrtum $\mathbb{L}_2$.}

 {\bf Satz 3}.\,\,{\it
 
 \noindent
 1. Das maximal Element des Spektrum $  \mathbb{L}_2$ ist $\frac{4}{\sqrt{5}}$,
 und
  $\hbox{\got k} (\alpha) =\frac{4}{\sqrt{5}}$ gilt dann und nur dann, wenn
 $\alpha$ und $ \frac{1+\sqrt{5}}{2}$
 sind \"aquivalent.

 \noindent
 2. Sei $\alpha$ irrational,
 und $\alpha, \frac{1+\sqrt{5}}{2}$
 seien miteinander \"aquivalent nicht. Gilt dann
 $\hbox{\got k}(\alpha) \le \frac{4}{\sqrt{17}}$,
 so
 ist
 $\left(
 \frac{4}{\sqrt{17}}
 ,
 \frac{4}{
 \sqrt{5}}
 \right)
 $
 eine L\"ucke im Spektrum $\mathbb{L}_2$.
 
  \noindent
 3. $\hbox{\got k} (\alpha) =\frac{4}{\sqrt{17}}$ gilt dann und nur dann, wenn
 $\alpha$ und $ \frac{1+\sqrt{17}}{2}$
 sind \"aquivalent.

 \noindent
  4. $\frac{4}{\sqrt{17}}$ ist ein
  isolierter Punkt der Menge
   $\mathbb{L}_2$.

  \noindent
    5. Das ganze Segment $\left[0,\frac{12}{21+\sqrt{15}}\right]$
        zum Spektrum $  \mathbb{L}_2$ geh\"ort.\footnote{Man kann leicht beweisen, dass ein gro\ss er Segment zum 
    $  \mathbb{L}_2$ geh\"ort.}

 }

 Beweis.
 
 F\"ur $\alpha = \frac{1+\sqrt{5}}{2} = [\overline{1}]$ ist
 $$
 \hbox{\got k}\left( \frac{1+\sqrt{5}}{2}\right) =
 \lim_{n\to \infty} \varkappa^4_n =
 \frac{4}{\sqrt{5}},
 $$
 und 
 f\"ur $\alpha = \frac{1+\sqrt{17}}{2} = [2;\overline{1,1,3}]$ ist
 $$
 \hbox{\got k} \left( \frac{1+\sqrt{17}}{2}\right) =
 \lim_{n\to \infty} \varkappa^4_{3n} =
 \lim_{n\to \infty} \varkappa^1_{3n}=
 \lim_{n\to \infty} \varkappa^2_{3n} =
 \frac{4}{\sqrt{17}}.
 $$ 
 Wir zeigen nun, dass aus $\alpha \not\sim \frac{1+\sqrt{5}}{2}$ folgt
 $
 \hbox{\got k} \left( \alpha\right) 
 \le
 \frac{4}{\sqrt{17}}.
 $
 Dann werden Behauptungen 1 und 2 bewiesen.

 Falls
 $\alpha \not\sim \frac{1+\sqrt{5}}{2},  \frac{1+\sqrt{17}}{2}$ 
  haben wir  die foldenden F\"alle 1 -  3 und 4.1 - 4.3.

 \noindent{\bf Fall 1$^0$.} Es gibt unendlich viele $n$ mit $ a_n \ge 5$. Dann f\"ur diese $n$ ist
 $\varkappa_n^4 \le \frac{4}{5} < \frac{4}{\sqrt{17}}.
 $
 
  \noindent{\bf Fall 2$^0$.} F\"ur alle hinreichend gro\ss en
  $n$ man hat $ a_n \le 4$, und
  es gibt unendlich viele $n$ mit $ a_n =4$. 
  Dann f\"ur diese $n$ ist
   $$
   \hbox{\got k}\left( \alpha\right) 
 \le
  \liminf_{n\to \infty, \, a_n = 4} \varkappa_n^4 \le \frac{4}{4+ 2\cdot [0;\overline{4,1}]} =
  \frac{4}{3+\sqrt{2}}< \frac{4}{\sqrt{17}}.
 $$

 \noindent{\bf Fall 3$^0$.}
  F\"ur alle hinreichend gro\ss en
  $n$ man hat $ a_n \le 4$, und
 es gibt unendlich viele $n$ mit $ a_n=2$. Dann ist
  $$
   \hbox{\got k} \left( \alpha\right) 
 \le
  \liminf_{n\to \infty, \, a_n = 2} \varkappa_n^1
  \le
  \frac{[1;\overline{1,4}|^2}{[2;\overline{1,4}]+[0;\overline{1,4}]}
  =
  \sqrt{2} - \frac{1}{2}
  < \frac{4}{\sqrt{17}}.
 $$

   \noindent{\bf Fall 4$^0$.}    F\"ur alle hinreichend gro\ss en
  $n$ man hat $ a_n \in \{1,3\}$

   \noindent{\bf Teilfall 4.1$^0$.} 
   Es gibt unendlich viele $n$ mit  $$a_{n-1} =a_{n} = 3. $$
   Dann ist
   $$
   \alpha_{n-1}^*\le [0;3,3] =\frac{13}{10}
   ,\,\,\,\,
   \alpha_n \le 4
   ,$$
   und
   $
      \varkappa_n^{1}
   \le
   \frac{39}{43}  < \frac{4}{\sqrt{17}}.$

   \noindent{\bf Teilfall 4.2$^0$.} 
    Es gibt unendlich viele  $n$ mit
$$
    a_{n-1}=3 , a_n =1, a_{n+1} = 3 .
 $$  
 Falls $a_{n-2} =3$  haben wir den Teilfall 4.1$^0$. 
 Sei $a_{n-2} = 1$. Dann ist 
$$
\varkappa_{n-1}^4 \le
\frac{4}{[3;1,3,3]+[0;1,1]} = \frac{104}{111}<\frac{4}{\sqrt{17}}.
$$

   \noindent{\bf Teilfall 4.3$^0$.} 
    Es gibt unendlich viele $n$ mit  
    $$
    a_{n-2}=a_{n-1}=1,a_n=3,a_{n+1}=a_{n+2}= a_{n+3} = 1.
 $$
 Dann ist 
$$
\varkappa_{n}^4 \le
\frac{4}{[3;1,1,1,1]+[0;1,1,3]} = \frac{140}{146}<\frac{4}{\sqrt{17}}.
$$

 Also in allen F\"allen  haben wir  
 $\hbox{\got k}  (\alpha) \le  \frac{140}{146},
 $ und damit sind  Behauptungen 1,2  und 3,4  bewiesen.

   Das Segment $\left[0,\frac{3}{21+\sqrt{15}}\right]$
  zum Lagrangeschen Spekrtum $\mathbb{L}$ geh\"ort.
  Wegen $ \frac{21+\sqrt{15}}{3} = [8;\overline{6,1}]+[0;\overline{6,1}]$,
  f\"ur alle $\lambda\in
  \left[0,\frac{3}{21+\sqrt{15}}\right]$ gibt es $\alpha$
  mit unendlich vielen 
 $a_n \ge 8$ und $\lambda (\alpha) = \lambda$.
  Nach Hilfssatz 13 (Behauptung 4)
  haben wir (\ref{lauf}). Sodass ist $\left[0,\frac{12}{21+\sqrt{15}}\right]\subset \mathbb{L}_2$.$\Box$

 {\bf 4.4. Das Spektrum $\mathbb{L}_2^*$.}

 {\bf Satz 4}.\,\,{\it
 
 \noindent
 1. Das maximal Element des Spektrum $  \mathbb{L}_2^*$ ist $\sqrt{5}$,
 und
  $\hbox{\got k}^* (\alpha) =\sqrt{5}$ gilt dann und nur dann, wenn
 $\alpha$ und $ \frac{1+\sqrt{5}}{2}$
 sind \"aquivalent.

 \noindent
 2. Sei $\alpha$ irrational,
 und $\alpha, \frac{1+\sqrt{5}}{2}$
 seien miteinander \"aquivalent nicht. Gilt dann
 $\hbox{\got k}^* (\alpha) \le \frac{3}{2}$,
 so
 ist
 $\left(
 \frac{3}{2}
 ,
 \sqrt{5}
 \right)
 $
 eine L\"ucke in dem Spektrum $\mathbb{L}_2^*$.
 
  \noindent
 3. F\"ur $ e=\sum_{k=0}^\infty \frac{1}{k!}$ ist $ \hbox{\got k}^* (e) =\frac{3}{2}$.

 \noindent
  4. $\frac{3}{2}$ ist ein H\"aufungspunkt der Menge
 $\mathbb{L}_2^*$.
 
  \noindent
  5. Das minimal Element des Spektrum $  \mathbb{L}_2^*$ ist $\frac{1}{2}$.

  \noindent
    6. Das ganze Segment $\left[\frac{1}{2},\frac{2}{3}\right]$
    zum Spektrum $  \mathbb{L}_2^*$ geh\"ort.

 }

 Beweis.

 Es ist klar, dass falls $\alpha \sim \frac{1+\sqrt{5}}{2}$,
 gilt
 $\hbox{\got k}^* (\alpha) =\sqrt{5}$.
 Betrachten wie eine irratiolale $\alpha\not\sim \frac{1+\sqrt{5}}{2}$.
 Wir finden unendlich viele $n$ mit $\hbox{\got l}^*_n (\alpha)\le \frac{3}{2}$.
  Es  gen\"ugt, die Beweise der Behauptungen 1 und 2 zu  beenden.

 F\"ur $\alpha \not\sim \frac{1+\sqrt{5}}{2}$ betrachten wir einige F\"alle.
 
 \noindent
 {\bf 
 Fall 1$^0$.}
 Es gibt unendlich viele $n $ mit  $a_n \ge 2$ und $a_{n+1}\ge 2$.
Dann
$
\varkappa_{n+1}^1 \le 1+\alpha_{n}^* \le \frac{3}{2}.
$

\noindent
 {\bf 
 Fall
 2$^0$.}
 F\"ur alle hinreichend gro\ss  en $n$,  entweder $a_n=1$ oder $a_{n+1}=1$. In diesem Fall
 der Kettenbruch f\"ur $\alpha$ ist
 \begin{equation}\label{ketten}
 \alpha =[a_0;\, ....\, ,\,\,
 \underbrace{1,...,1}_{r_1}, a_{n_1},
  \underbrace{1,...,1}_{r_2},a_{n_2},
  \underbrace{1,...,1}_{r_3},
  \,\,\,\,
  ...
  \,
  ,
  \,\,\,\,
  \underbrace{1,...,1}_{r_j}, a_{n_j},
  \underbrace{1,...,1}_{r_{j+1}}, ...
 ],\,\,\,\, r_j \ge 1,\,\,\, a_{n_j} \ge 2.
 \end{equation}

 \noindent
 {\bf 
 Teilfall
 2.1$^0$.} Es gibt unendlich viele $j$ mit $ r_j =1$.
 Dann gibt es unendlich viele $n$ mit
 $$
 a_{n-1} \ge 2, a_n = 1, a_{n+1} \ge 2.
 $$

  \noindent
 {\bf 
 Teilfall
 2.1.1$^0$.} $a_{n-1} \le 6$ oder $a_{n+1} \le 6$. Dann aus (\ref{man2}) folgt
  $\min (\varkappa^1_{n-1}, \varkappa^1_{n+1})  \le \frac{3}{2}$.

 \noindent
  {\bf 
 Teilfall
 2.1.2$^0$.} $\min (a_{n-1} , a_{n+1} ) \ge 7$. Dann 
 gelten $ \alpha_{n} = [1; a_{n+1},...] \le \frac{8}{7}$ und
 $\alpha_{n-1}^* = [0;a_{n-1},...] \le \frac{1}{7}$,
 und 
 aus der Definition (\ref{kk3})  haben wir
 $\varkappa_n^3 \le \frac{9}{7}< \frac{3}{2}.$
 
 \noindent
 {\bf 
 Teilfall
 2.2$^0$.} 
  F\"ur alle hinreichend gro\ss  en $j$,
  man hat
 $r_j\ge 2$.

 \noindent
 {\bf 
 Teilfall
 2.2.1$^0$.} 
 Es gibt unendlich viele $j$ mit $ r_j = 2$.
 Dann gibt es unendlich viele $n$ mit
 $$
 a_{n-2} =a_{n-1}=1, a_n=a \ge 2, a_{n+1} =a_{n+2} = 1, a_{n+3}=b \ge 2,
 a_{n+4}=a_{n+5}= 1. 
 $$
 
 Falls $ a\ge b$,  haben wir
 $\alpha_{n+3} = [b;1,1,...]\le [b;1,1,1] = b+\frac{2}{3}$
 und
 $\alpha_{n+2}^*=[0; 1,1,a,1,1,...]\le[0;1,1,b,1,1] 
=\frac{2b+3}{4b+4}$. Dann gilt
$\varkappa^1_{n+3} \le \frac{18b^2+ 15b-7}{12b^2+26b +17}< \frac{3}{2}$
f\"ur alle $b$.
 
 Falls $ a\ge b$, beweisen wir (siehe Bemerkung 3), dass $\varkappa^2_n < \frac{3}{2}$
f\"ur alle $a$.

 \noindent
 {\bf 
 Teilfall
 2.2.2$^0$.}   F\"ur alle hinreichend gro\ss  en $j$,
  man hat
 $r_j\ge 3$. 
 
 \noindent
 {\bf 
 Teilfall
 2.2.2.1$^0$.}
 Es gibt unendlich viele $n$ mit
 $$ a_{n-3} = a_{n-2}=a_{n-1} =1, a_n = a, a_{n+1} =a_{n+2}=a_{n+3}=1
 $$
 und
 $2\le  a\le 15$.
 In diesem Falle  haben wir
 $$
 \alpha_n =[a;1,1,1,...]\le [a;,1,1,1]= a+\frac{2}{3} \le \frac{47}{3},\,\,\,\,
 \alpha_{n-1}^* =[0;1,1,1,...] \le \frac{2}{3},
$$
 und
 $
 \varkappa_n^1 \le \frac{220}{147}< \frac{3}{2}.
 $

 \noindent
 {\bf 
 Teilfall
 2.2.2.2$^0$.}
  Es gibt unendlich viele $n$ mit
   $$ a_{n-3} = a_{n-2}=a_{n-1} =1, a_n = a, a_{n+1} =a_{n+2}=a_{n+3}=1,
   a_{n+4}=b,
   a_{n+5} =a_{n+6}=a_{n+6}=1,
   $$
   und $ a\ge 16, b\ge 16$.
 In diesem Falle  haben wir
 $$
 \alpha_{n+3}=[1;b,1,1,1,...]\le [1;16,1,1,1,1] =\frac{88}{83},\,\,\,\,
 \alpha_{n+2}^* = [0;1,1,a,1,1,1,...]\le [0;1,1,16,1,1,1,1] = \frac{88}{171},
 $$
 und 
 $\varkappa_{n+3}^3 \le 
 \frac{\left(2\cdot \frac{88}{171}+1\right)\left(2\cdot\frac{88}{83}-1\right)}{\frac{88}{171}+\frac{88}{83}}< \frac{3}{2}$.

 \noindent
 {\bf 
 Teilfall
 2.2.2.3$^0$.}
 Es gibt unendlich viele $n$ mit
  $$a_{n-4} = a_{n-3} = a_{n-2}=a_{n-1} =1, a_n = a, a_{n+1} =a_{n+2}=a_{n+3}=a_{n+4}=1.
 $$

 Falls $a\le 19$,
  haben wir
 $$
 \alpha_n =[a;,1,1,1,1,...]\le [a;1,1,1,1,1] = a+\frac{5}{8}\le \frac{157}{8}
 ,\,\,\,\,
 \alpha_{n-1}^* = [0;1,1,1,1,...]\le [0;1,1,1,1,1] =\frac{5}{8},
 $$
 sodass
 $\varkappa^1_n \le \frac{13\cdot149}{8\cdot 163}< \frac{3}{2}$.

 Falls $a\ge 20$,
  haben wir
 
 $$
 \alpha_{n-1} =[1;a,1,1,1,1,...]\le 1+\frac{1}{a+\frac{3}{5}}\le \frac{108}{103}
 ,\,\,\,\,
 \alpha_{n-2}^* = [0;1,1,1,...]\le [0;1,1,1] =\frac{2}{3},
 $$
 sodass
 $\varkappa^3_{n-1} \le \frac{7\cdot113}{530}< \frac{3}{2}$.

 Also in allen F\"allen  haben wir unendlich viele $n$
 mit $\hbox{\got l}_n^* \le \frac{3}{2}$
 gefunden, und damit sind die Behauptungen 1,2  bewiesen.

 Sei $\{a_j\}_{j=0}^\infty$ eine ganzzahlige Sequenz mit $\lim_{j\to +\infty} a_j= +\infty$.
 Betrachten wir $ \alpha = [a_0;a_1,...,a_t, \overline{a_j,1,1}_{j=t+1}^\infty]$.
  Es ist klar aus dem Beweis des Teilfalles 2.2.1$^0$, dass
  $\lim_{j\to \infty}\varkappa^1_{t+1+3j} =\frac{3}{2}$ und $\hbox{\got k}^*(\alpha) = \frac{3}{2}$.
  Als $e = [2;1,\overline{2+2j,1,1}_{j=0}^\infty]$, gilt auch $\hbox{\got k}^*(e) = \frac{3}{2}$,
  und die Behauptung 3 ist bewiesen.

  Behauptung 4 ist klar.
  
  F\"ur jede irrationale $\alpha$ 
  aus (\ref{man1},\ref{man2})
  folgt $\hbox{\got k} (\alpha) \ge \frac{1}{2}$.
  Sei $\{b_j\}_{j=0}^\infty$ eine ganzzahlige Sequenz mit $\lim_{j\to +\infty} b_j= +\infty$.
 Definiren wir $\beta =[0;\overline{b_j,2}_{j=0}^\infty]$.
 Dann gilt $\lim_{j\to \infty} \varkappa^1_{2j} =\frac{1}{2}$. Sodass
 $\hbox{\got k}^*(\beta ) = \frac{1}{2} =\min \, \mathbb{L}_2^*$ und wir haben Behauptung 5 bewiesen.

 Nun  sollen wir die Behauptung 6 beweisen. 
 Hier werwenden  wir eine Methode von M. Hall \cite{Hall}. Wir nutzen eine Verallgemeinerung aus \cite{alex,astels,mshch}.

Bezeichnen wir mit ${\cal F}$ die Menge aller endlichen und unendlichen Kettenbr\"uche
$\xi =[0;a_1,...,a_n,...]$ aus dem Segment $\left[0,\frac{1}{3}\right]$,
deren Teilnenner nicht gleich zu 2 sind:
 $$
 {\cal F}
=
\{\xi\in  
\left[0,\frac{1}{3}\right]:\,\,
\xi \,\,\,\text{kann durch einen endlichen oder unendlichen  Kettenbruch}
$$
$$
 \,\,\,\xi =[0;a_1,...,a_n,...]\,\,\,
  \text{dargestellt werden mit}\,\,\, a_k \neq 2 \,\, \forall \,\, k\}.
 $$
 Dann
 $$
  {\cal F}
= 
 \left[0,\frac{1}{3}\right]
 \setminus \left(
 \bigcup_{\nu=1}^\infty
 \Delta_\nu
 \right),
 $$
 wo $\Delta_\nu$ sind  offene Intervalle.
 Diese Intervalle k\"onnen so geordnet werden, dass das Folgendes gilt
 (siehe Hilfssatz 4.3 aus \cite{astels} mit $b_1,b_2,b_3,b_4,...
= 1,3,4,5,...$).
 F\"ur jede $T$ man hat
 $$
  \left[0,\frac{1}{3}\right]
 \setminus \left(
 \bigcup_{\nu=1}^T
 \Delta_\nu
 \right)=
 \bigcup_{j=1}^{T+1} I_{T,j}.
 $$
 wo  $I_{T,j}$ abgeschlossene Segmente sind, und falls
 $$
 \Delta_{T+1} \subset I_{T,j^*},
 $$
 gilt
 $$
 I_{T,j^*} = M_1\sqcup \Delta_{T+1}\sqcup M_2,
 $$
 wo f\"ur die 
 L\"angen der Intervalle $M_k , k =1,2$ man hat
 $$
 \min ( \text{ L\"ange}\, M_1,  \text{ L\"ange}\, M_2) \ge \tau \cdot  \text{ L\"ange}\, \Delta_{T+1},\,\,\,\,
 \tau = 2.
 $$
 
 Definiren wir
 $$
 H(x,y) =\frac{(1+x)(1+y)}{2+x+y}.
 $$
 Dann
 $$
 \frac{\partial H /\partial x}{\partial H/\partial y} = \left( \frac{1+y}{1+x}\right)^2,\,\,\,
 \frac{\partial H /\partial y}{\partial H/\partial x} = \left( \frac{1+x}{1+y}\right)^2,
 $$
 und
 $$
 \max_{0\le x,y\le \frac{1}{3}} \,
 \max\,\left(
 \left|\frac{\partial H /\partial x}{\partial H/\partial y}\right|,
 \left|\frac{\partial H /\partial y}{\partial H/\partial x}\right| 
 \right)
 \le \frac{16}{9} < 2 = \tau.
 $$
 Nun nach Hilfssatz 2 aus \cite{alex} folgt, dass das Bild der Menge ${\cal F}\times{\cal F}$ 
 unter $H$ ein  abgeschlossenes Intervall ist:
 $$
 H({\cal F},{\cal F}) = \left[H(0,0), H\left(\frac{1}{3},\frac{1}{3}\right)\right] 
 =
 \left[\frac{1}{2},\frac{2}{3}\right] 
 .
 $$
 Bemerken wir, dass falls $a_n = 2$, man hat
 $$
 \varkappa_n^1 =
 \varkappa_n^2 =
 H(\alpha_{n-1}^*, 1/\alpha_{n+1}),
 $$
 und falls $a_n \neq 2$, aus (\ref{man1},\ref{man2}) folgen 
 $$
 \min (\varkappa_n^1, \varkappa_n^2,\varkappa_n^3)\ge \frac{2}{3}.
 $$
 Nun  nehmen
wir
 $k \in \left[\frac{1}{2},\frac{2}{3}\right] $.
 Betrachten wir
 zwei irrationale Zahlen
 $$
 x =[0;x_1,...,x_n,...],
  y =[0;y_1,...,y_n,...] \in {\cal F}
 $$
 mit
 $$
 H(x,y) = k
 $$
 und den Kettenbruch
 $$
 \alpha =
 [0; \underbrace{x_1,2,y_1}, \,\,\underbrace{x_2,x_1 2,y_1,y_2},\,\,
 \underbrace{x_3,x_2,x_1 2,y_1,y_2,z_3},\,\,\,
 ...\,\,\, , \underbrace{x_k, x_{k-1}, , ..., x_2,x_1,2,y_1,y_2,  ...  ,y_{k-1},y_k}, \,\,\,...
 ]
 .
 $$ 
 Wegen $x_n, y_n \neq 2$,
 es ist klar, dass $\hbox{\got k}^* (\alpha ) = k$. $\Box$

 {\bf Bemerkung 6.}
 In \cite{alex} Hilfssatz 2 wurde ohne Beweis formuliert.
 Aber der Beweis folgt sofort aus dem Argument aus \cite{mshch}.

Der Autor dankt 
 Oleg N. German, Igor D. Kan und Igor P. Rochev f\"ur die fruchbaren Diskussionen.

 \end{document}